# On the asymptotic Plateau problem in Cartan–Hadamard manifolds.

Graham Smith*

**Abstract:** We prove dynamical stability of a natural class of hypersurface laminations defined over Cartan–Hadamard manifolds of pinched curvature. We achieve this by providing a complete solution to the asymptotic Plateau problem for immersed surfaces of constant extrinsic curvature in Cartan–Hadamard manifolds, proposed by Labourie in [25], together with its natural higher-dimensional generalisation.

## 1 - Introduction.

**1.1 - The laminated Plateau problem.** Let $X := X^{d+1}$ be a $(d+1)$-dimensional Cartan–Hadamard manifold of sectional curvature bounded above by $-1$, let $SX$ denote its unit sphere bundle, and let $\partial_\infty X$ denote its ideal boundary. When $d = 1$, it is well known that Cheeger's homeomorphism (see [16] and Figure 1.1.1) naturally identifies $SX$ with the set of ordered triples of distinct points in $\partial_\infty X$, yielding the following useful fact: any homeomorphism $\alpha : \partial_\infty X \to \partial_\infty X'$ between the ideal boundaries of two such Cartan–Hadamard surfaces induces a natural homeomorphism $A : SX \to SX'$ of their unit circle bundles which conjugates their geodesic foliations. We see here a simple instance of a general principal: much of the geometry and dynamics of Cartan–Hadamard manifolds is encoded in their ideal boundaries. We present an analogue of this result for higher-dimensional laminations in $SX$.

 We first describe the hypersurfaces of which these laminations are formed. Let $(Y, e)$ be an immersed hypersurface in $X$, let $\hat{e}$ denote its unit normal vector field compatible with its orientation, and let $\text{I}_e$, $\text{II}_e$ and $\text{III}_e$ denote its three fundamental forms. We say that $(Y, e)$ is **quasicomplete** whenever it is complete with respect to the metric $\text{I}_e + \text{III}_e$, and we say that it is **infinitesimally strictly convex** (ISC) whenever $\text{II}_e$ is everywhere positive-definite.

 In [32] and [35], we introduced the concept of special Lagrangian curvature (see Section 3.1) for ISC hypersurfaces. This curvature notion arises naturally from the contact structure of $SX$ and serves as a higher-dimensional analogue of 2-dimensional extrinsic curvature, that is, the determinant of the shape operator. The reader may consult the recent work [19] of Harvey–Lawson for a discussion thereof within the broader context of the special lagrangian potential equation.

 For $0 < k < 1$, we define a $k$-**hypersurface** to be a quasicomplete, ISC, immersed hypersurface $(Y, e)$ of constant special Lagrangian curvature equal to $k$. When $d = 2$, $k$-hypersurfaces are simply surfaces of constant extrinsic curvature equal to $k$. These are Labourie's $k$-surfaces, whose properties he studied in [24], [25] and [26], obtaining a variety of deep results.

 The laminations of $SX$ of interest to us are constructed out of 1-jets of $k$-hypersurfaces and complete geodesics as follows. For $0 < k < 1$, let $\mathcal{S}_k(X)$ denote the set of triplets of the form $(Y, \hat{e}, p)$, for some $k$-hypersurface $(Y, e)$, and some point $p \in Y$, where we recall from above that $\hat{e}$ denotes the unit normal vector field over $e$. Let $\mathcal{T}$ denote the set of the triplets of the form $(Y, \hat{e}, p)$, where $(Y, \hat{e})$ is a covering of the unit, normal sphere bundle of some complete geodesic in $X$ and $p \in Y$. For all $0 < k < 1$, we denote

$$\overline{\mathcal{S}}_k(X) := \mathcal{S}_k(X) \cup \mathcal{T}(X), \tag{1.1}$$

we identify elements which are equivalent up to reparametrisation, and we furnish this space with the Cheeger–Gromov topology (see Section 4.1). $\overline{\mathcal{S}}_k(X)$ carries a natural laminated structure, outlined in the 2-dimensional case in [25]. We view $\mathcal{T}(X)$ as a lifting of the geodesic flow of $SX$. It is a finite-dimensional sublamination of $\overline{\mathcal{S}}_k(X)$, and its interpretation as the topological boundary of $\mathcal{S}_k(X)$ is justified by Theorem 4.1.2.

---

* Universidade Federal do Rio de Janerio, Rio de Janeiro, Brazil

On the asymptotic Plateau problem in Cartan-Hadamard manifolds.

We prove dynamical stability of the lamination $\overline{\mathcal{S}}_k(X)$. For all $0 < k < 1$, we define $\Phi : \overline{\mathcal{S}}_k(X) \to \partial_\infty X$ by

$$\Phi(Y, \hat{e}, p) := (\text{Hor} \circ \hat{e})(p), \tag{1.2}$$

where $\text{Hor} : SX \to \partial_\infty X$ denotes the **horizon map**, that is, for all $\xi \in SX$,

$$\text{Hor}(\xi) := \lim_{t \to +\infty} \gamma_\xi(t), \tag{1.3}$$

where $\gamma_\xi : \mathbb{R} \to X$ denotes the unique geodesic such that $\dot{\gamma}_\xi(0) = \xi$.

**Theorem 1.1.1, Dynamical stability**

*Let $X$ and $X'$ be $(d+1)$-dimensional Cartan–Hadamard manifolds of sectional curvature pinched between $-K$ and $-1$ for some $K \geqslant 1$. For all $0 < k, k' < 1$, and for every homeomorphism $\psi : \partial_\infty X \to \partial_\infty X'$, there exists a unique homeomorphism $\Psi : \overline{\mathcal{S}}_k(X) \to \overline{\mathcal{S}}_{k'}(X')$, which sends leaves to leaves, such that the following diagram commutes.*

$$\begin{array}{ccc} \overline{\mathcal{S}}_k(X) & \xrightarrow{\Psi} & \overline{\mathcal{S}}_{k'}(X') \\ \Phi \downarrow & & \downarrow \Phi \\ \partial_\infty X & \xrightarrow{\psi} & \partial_\infty X' \end{array} \tag{1.4}$$

*Furthermore, $\Psi$ maps $\mathcal{T}(X)$ into $\mathcal{T}(X')$.*

**Remark 1.1.1.** Theorem 1.1.1 is proven in Section 4.6.

Theorem 1.1.1 should be understood in the context of foliated Plateau problems, introduced by Gromov in [14] and [15]. To paraphrase Gromov, by viewing constant curvature submanifolds "not as individuals, but as members of a community," techniques of dynamical systems theory may be brought to bear on their study. In recent years, striking applications of dynamical systems theory to the study of minimal surfaces have indeed been developed, most notably by Coda–Neves' in their proof [10] of the Wilmore conjecture, by Coda–Neves–Song in their work [11] on the equidistribution of minimal surfaces in compact 3-dimensional manifolds, and by Calegari–Marques–Neves in their study [9] of the asymptotic counting of compact minimal surfaces in negatively curved 3-dimensional manifolds.

From a more dynamical perspective, in a remarkable series of papers [24], [25] and [26], Labourie fruitfully applied Gromov's vision to the study of constant extrinsic curvature surfaces in Cartan–Hadamard manifolds. In addition, in [27], he used concepts analogous to those developed in these papers to highlight the equidistribution properties underlying Calegari–Marques–Neves aforementioned work (we revisit these ideas in greater detail in [1]).

We view Theorem 1.1.1 as a development of Labourie's programme. Indeed, in the aforementioned papers, Labourie addresses the case of surface laminations in 3-dimensional Cartan–Hadamard manifolds which are geometrically finite in a certain strong sense. In [25] he shows that, in this case, $\overline{\mathcal{S}}_k(X)$ possesses the hyperbolic properties of the geodesic flow. In particular, it is stable under *small* perturbations of the ambient space. Theorem 1.1.1 extends Labourie's dynamical stability result to one which is canonical, is valid for large perturbations, and holds over a more general class of ambient spaces, and in arbitrary dimension.

**1.2 - The asymptotic Plateau problem.** Theorem 1.1.1 follows from the complete solution of an asymptotic Plateau problem for $k$-hypersurfaces in Cartan–Hadamard manifolds of pinched curvature, formulated in the case of surfaces by Labourie in [24]. The asymptotic Plateau problem for constant extrinsic curvature hypersurfaces in hyperbolic space was first addressed by Rosenberg–Spruck in [30], and extended to more general curvature functions by Guan–Spruck in, for example, [17]. In these works, the hypersurface is prescribed by its boundary curve, and the question of well-posedness can be addressed by prescribing the orientation of this curve. In [25], Labourie noted that well-posedness is also guaranteed by instead prescribing

On the asymptotic Plateau problem in Cartan-Hadamard manifolds.

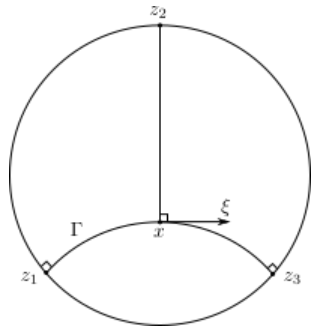

**Figure 1.1.1 - The Cheeger homeomorphism -** This homeomorphism identifies the unit circle bundle over $X$ with an open subset of $\partial_\infty X^3$ as follows. The unit vector $\xi$ over the point $x$ identifies with the two end-points $z_1$ and $z_3$ of the geodesic $\Gamma$ that it defines, together with the unique point $z_2$ lying to the left of $\Gamma$ which projects onto $x$.

one of the two connected components of the complement of this curve in the ideal boundary. This naturally led to a novel formulation of the asymptotic Plateau problem in terms of the asymptotic Gauss map, which has since shown itself to better reflect the properties of $k$-surfaces, and the ideas developed in this pioneering work have yielded deep results across a variety of fields (see, for example, [4], [5], [6], [22], [23], [26] and [31]).

Before stating Labourie's asymptotic Plateau problem, it is worth noting that the concept of smooth immersion $e : Y \to X$ does not require a differential structure over the domain $Y$ in order to be defined. Indeed, given a topological manifold $Y$, we say that a function $e : Y \to X$ is a **smooth immersion** whenever there exists a smooth atlas over $Y$ with respect to which $e$ has this property. Trivially, any two such atlases are contained within the same maximal atlas, which we call the atlas **induced** by $e$. We define an **immersed hypersurface** in $X$ to be a pair $(Y, e)$, where $Y$ is a $d$-dimensional topological manifold and $e : Y \to X$ is a smooth immersion. When working with immersed hypersurfaces, we will always use the induced atlas without comment. In addition, we will always assume that the hypersurfaces we work with are oriented.

Given a $k$-hypersurface $(Y, e)$, we define its **asymptotic Gauss map** by

$$\phi_e := \text{Hor} \circ \hat{e}, \tag{1.5}$$

where Hor denotes the horizon map defined above. By standard properties of convex subsets of Cartan–Hadamard manifolds (see, for example, [3]), this function is a local homeomorphism. Labourie's asymptotic Plateau problem now asks whether $k$-hypersurfaces in Cartan–Hadamard manifolds can be uniquely prescribed by their asymptotic Gauss maps. More precisely, following [25], we make the following definition.

**Definition 1.2.1**

*Let $X$ be a $(d+1)$-dimensional Cartan–Hadamard manifold and let $\partial_\infty X$ denote its ideal boundary. We define an **asymptotic Plateau problem** in $\partial_\infty X$ to be a pair $(Y, \phi)$, where $Y$ is a connected, oriented topological manifold, and $\phi : Y \to \partial_\infty X$ is a local homeomorphism. For all $k > 0$, we say that a $k$-hypersurface $(Y, e)$ **solves** the asymptotic Plateau problem $(Y, \phi)$ whenever its asymptotic Gauss map $\phi_e$ satisfies*

$$\phi_e = \phi. \tag{1.6}$$

**Remark 1.2.1.** The classical Plateau problem at infinity is shown to be a special case of this problem as follows. Let $\Gamma$ be a Jordan sphere in $\partial_\infty X$, that is, a subset homeomorphic to a $(d-1)$-dimensional sphere. Let $Y$ be one of the connected components of its complement, and let $\phi : Y \to \partial_\infty X$ denote the canonical embedding. If, for $k > 0$, $e : Y \to X$ solves $(Y, \phi)$ then, reasoning as in Lemmas 3.4.1 and 3.4.2, below, we show that $e$ is in fact a complete embedding with ideal boundary $\Gamma$, and thus solves the classical Plateau problem at infinity defined by this curve.

In [25], Labourie obtains partial existence, non-existence and uniqueness results for solutions of the 2-dimensional asymptotic Plateau problem. The general case, however, even in 2-dimensions, has since remained open. On the one hand, Labourie's results require the ambient space to be of bounded geometry, which is a restrictive condition on its structure at infinity. On the other hand, various straightforward Plateau problems lie beyond the scope of Labourie's work. For example, the case where $X = \mathbb{H}^3$, $Y = \mathbb{C} \setminus \{0\}$, and

$$\phi : Y \to \hat{\mathbb{C}}; z \mapsto e^z \tag{1.7}$$

is not addressed by [25], nor are problems constructed using more complicated entire functions, nor topologically similar problems in more general spaces.

Our second main result provides a complete solution to Labourie's asymptotic Plateau problem, over a broader class of ambient manifolds, and in arbitrary dimension. We require one final definition in order to address exceptional cases. For any non-negative integer $m$, we will say that the asymptotic Plateau problem $(Y, \phi)$ is of **finite type** $m$ whenever $\phi$ is a topological cover of $\partial_\infty X \setminus F$, for some finite subset $F$ of cardinality $m$.

**Theorem 1.2.2, Existence, non-existence and uniqueness**

Let $X$ be a $(d+1)$-dimensional Cartan–Hadamard manifold of sectional curvature pinched between $-K$ and $-1$, for some $K \geqslant 1$. Let $(Y, \phi)$ be an asymptotic Plateau problem in $\partial_\infty X$.

(1) If $(Y, \phi)$ is of finite type 0, 1 or 2 then, for all $0 < k < 1$, there exists no $k$-hypersurface solving $(Y, \phi)$.

(2) If $(Y, \phi)$ is not of finite type 0, 1 or 2 then, for all $0 < k < 1$, there exists a unique $k$-hypersurface solving $(Y, \phi)$.

**Remark 1.2.2.** Theorem 1.2.2 is proven in Section 4.5.

**1.3 - Summary and open problems.** This paper is structured as follows.

**Chapter 2:** We introduce the concept of developed Cartan–Hadamard ends in $X$. We define their ideal boundaries and show that they yield asymptotic Plateau problems in $\partial_\infty X$. We associate to every asymptotic Plateau problem $(Y, \phi)$ in $\partial_\infty X$ a certain continuous section $\omega_\phi$ of a certain topological bundle $\phi^*\mathrm{HB}(X)$, which generalises to the non-constant curvature case the metric described by Kulkarni–Pinkall in [20], and which we thus call the Kulkarni–Pinkall section. This section encodes in a local manner, not only the global geometry of $(Y, \phi)$, but also that of every developed Cartan–Hadamard end in $X$ having $(Y, \phi)$ as ideal boundary. Using this section we obtain in Theorem 2.8.2 an a priori $C^0$ estimate for infinitesimally strictly convex, quasicomplete immersed hypersurfaces in $X$.

**Chapter 3:** We recall the definition and basic properties of special Lagrangian curvature. We estimate the special lagrangian angles of the level sets of the height functions of certain Cartan–Hadamard ends. Using a suitably adapted version of Omori's maximum principle (see Appendix B), we prove in Lemma 3.3.2 a partial uniqueness result. We conclude by establishing criteria for one developed Cartan–Hadamard end to be contained within another.

**Chapter 4:** We recall the precompactness theory of $k$-surfaces described by Labourie in [24] and extended to the case of $k$-hypersurfaces by the author in [35]. Using this theory together with the a priori estimate of Theorem 2.8.2, we obtain in Theorem 4.2.1 a monotone convergence result for families of $k$-hypersurfaces. We recall the theory of non-linear Plateau problems and convex cobordisms developed by the author in [36]. In Section 4.5, we combine this theory with Theorem 4.2.1 and Lemma 3.3.2 to yield Theorem 1.2.2. Finally, Theorem 1.1.1 is proven in Section 4.6.

Finally, we summarise a number of interesting byproducts of our work and open problems that remain.

(1) Theorem 1.2.2 may also be expressed in terms of the existence and uniqueness of complete special legendrian immersions of the unit sphere bundle (see Appendix C). In particular, by describing the special legendrian geometry of this bundle in terms of the ideal boundary, Theorem 1.1.1 may be viewed as a manifestation of the holographic principle.

(2) Our techniques yield a generalisation of the construction [20] of Kulkarni–Pinkall to the non-constant curvature case (see Appendix D). In particular, one wonders whether, like in the constant curvature case, the Cartan–Hadamard ends constructed in Theorem D.1 are maximal (see Remark D.2).

(3) We expect Theorem 1.2.2 to continue to hold also for surfaces of constant *intrinsic* curvature equal to $\kappa$, for all $-1 < \kappa < 0$. In this manner, we would obtain laminations of Cartan–Hadamard 3-manifolds by hyperbolic surfaces. It is natural to ask under which hypotheses these surfaces would be complete, as opposed to merely quasicomplete.

(4) In Section 1.2 of [25], Labourie speculates on a number of dynamical properties of the space of $k$-surfaces, many of which have to date remained unstudied and present fascinating open problems. In particular, in the light of the work [9] of Calegari–Marques–Neves and our own work [1] in collaboration with S. Alvarez and B. Lowe, it is natural to ask whether the entropy proposed by Labourie, or a suitable modification thereof, can be estimated, whether this estimate is sharp, and whether rigidity holds in the case of equality. It is also natural to ask to what extent these dynamical properties can also be studied in the higher-dimensional case.

**1.4 - Acknowledgements.** The author is grateful to François Labourie for having initially brought this problem to his attention. The author is also grateful to Sébastien Alvarez, François Fillastre, Martin Kilian, and Andrea Seppi for helpful comments made to earlier drafts of this paper.

## 2 - Asymptotic Plateau problems and Cartan–Hadamard ends.

**2.1 - Preliminary constructions in Cartan–Hadamard manifolds.** We begin by recalling some elementary definitions of the theory of Cartan–Hadamard manifolds. Let $X := X^{d+1}$ be a $(d+1)$-dimensional Cartan-Hadamard manifold, let $TX$ denote its tangent bundle and let $SX \subseteq TX$ denote its unit sphere bundle. Let $\partial_\infty X$ denote the **ideal boundary** of $X$, which we recall is defined to be the set of equivalence classes of complete, unit speed geodesic rays, where two such rays are deemed equivalent whenever they remain within bounded distance of one another. The union $X \cup \partial_\infty X$ is furnished with the **cone topology**, which we recall is defined as follows (c.f. [3]). For all $x \in X$, for all $\xi \in S_x X$, for all $\theta \in ]0, 2\pi[$, and for all $r > 0$, the **open truncated cone** of angle $\theta$ and inner radius $r$ about $\xi$ is defined by

$$\mathrm{C}(\xi, \theta, r) := \{\mathrm{Exp}_x(t\mu) \mid \mu \in S_x X, \ t > r, \ \langle \mu, \xi \rangle > \cos(\theta)\}. \tag{2.1}$$

For all such $\xi$, $\theta$ and $r$, the **ideal boundary** $\partial_\infty \mathrm{C}(\xi, \theta, r)$ of $\mathrm{C}(\xi, \theta, r)$ is defined to be the set of equivalence classes of rays which eventually lie in this cone. The collection of all sets of the form

$$\mathrm{C}(\xi, \theta, r) \cup \partial_\infty \mathrm{C}(\xi, \theta, r)$$

together with the topology of $X$ then defines a basis of a compact Hausdorff topology over $X \cup \partial_\infty X$ which is homeomorphic to that of the closed unit ball in $\mathbb{R}^{d+1}$. In particular, for any complete, unit speed geodesic ray $\gamma : [0, \infty[ \to X$ with equivalence class $[\gamma]$, we write

$$\underset{t \to \infty}{\mathrm{Lim}}\, \gamma(t) = [\gamma]. \tag{2.2}$$

We will use this notation frequently throughout the sequel.

For all $\xi \in SX$, we define $\gamma_\xi : \mathbb{R} \to X$ to be the unique geodesic such that

$$\dot\gamma_\xi(0) := \xi. \tag{2.3}$$

We then define the **horizon map** $\mathrm{Hor} : SX \to \partial_\infty X$ by

$$\mathrm{Hor}(\xi) := \underset{t \to \infty}{\mathrm{Lim}}\, \gamma_\xi(t). \tag{2.4}$$

This function is continuous and, for any given point $x$ of $X$, since every class of $\partial_\infty X$ contains a unique geodesic ray emanating from $x$, it restricts to a homeomorphism

$$\mathrm{Hor}_x : S_x X \to \partial_\infty X.$$

It is well known that $\partial_\infty X$ does not in general carry a natural differential structure, since the composition $\mathrm{Hor}_y^{-1} \circ \mathrm{Hor}_x$ is typically not smooth. However, given any choice of base point $x$, there trivially exists a unique differential structure over $\partial_\infty X$ making $\text{Hor}_x$ into a diffeomorphism. This observation will prove useful in the sequel.

We conclude this section by defining families of special subsets of $X$ which we will use frequently throughout the rest of the paper. For all $\xi \in SX$, we define the **open half-space** centred on $\xi$ by

$$\text{HS}(\xi) := C(\xi, \pi/2, 0) = \{\text{Exp}_x(\nu) \mid \nu \in T_xX,\ \langle \nu, \xi \rangle > 0\}, \tag{2.5}$$

and the **asymptotic ball** centred on $\xi$ by

$$\text{AB}(\xi) := \partial_\infty \text{HS}(\xi) = \{\text{Hor}_x(\nu) \mid \nu \in S_xX,\ \langle \nu, \xi \rangle > 0\}. \tag{2.6}$$

In addition, for all $\xi \in SX$ and for all $t \in \mathbb{R}$, we denote

$$\text{HS}_t(\xi) := \text{HS}(\dot{\gamma}_\xi(t)), \tag{2.7}$$

where $\gamma_\xi$ is again as in (2.3). Finally, for all $\xi \in SX$, we denote by $\text{CE}(\xi)$ the complement of the convex hull in $X$ of the complement of $\text{AB}(\xi)$. When the sectional curvature of $X$ is pinched between $-K$ and $-1$, it follows from the work [2] of Anderson that, for all such $\xi$,

$$\partial_\infty \text{CE}(\xi) = \text{AB}(\xi). \tag{2.8}$$

We call $\text{CE}(\xi)$ the **canonical end** centred on $\xi$ for reasons that should become clear presently.

**Lemma 2.1.1**

Let $X$ be a Cartan–Hadamard manifold with sectional curvature pinched between $-K$ and $-1$. There exists $t_0 > 0$ which only depends on $K$ such that, for all $\xi \in SX$,

$$HS_{t_0}(\xi) \subseteq CE(\xi) \subseteq HS(\xi). \tag{2.9}$$

**Proof:** First observe that $\partial \text{HS}(\xi)$ consists of geodesics passing through $x$ so that, by convexity, $x \in \text{CE}(\xi)^c$ and, by convexity again, $\text{HS}(\xi)^c \subseteq \text{CE}(\xi)^c$, from which the second inclusion follows. To prove the first inclusion, let $\gamma_\xi$ be as in (2.3). The construction of Section 2 of [2] yields a convex subset $Y$ of $X$ such that

$$\text{AB}(\xi)^c \subseteq \partial_\infty Y.$$

Furthermore, there exists $t_0 > 0$, which only depends on $K$, with the property that $Y$ can be chosen in such a manner that $\gamma(t_0) \in \partial Y$ and, near $\gamma(t_0)$, $\partial Y$ coincides with the geodesic sphere $\partial B_{t_0}(x)$. It follows by convexity that

$$Y \subseteq \text{HS}_{t_0}(\xi)^c,$$

so that

$$\text{CE}(\xi)^c \subseteq Y \subseteq \text{HS}_{t_0}(\xi)^c,$$

and the first inclusion follows. $\square$

**2.2 - Asymptotic Plateau problems.** We now recall the definition and basic properties of asymptotic Plateau problems. Let $X$ be a Cartan–Hadamard manifold with ideal boundary $\partial_\infty X$. We define an **asymptotic Plateau problem** in $\partial_\infty X$ to be a pair $(Y, \phi)$, where $Y := Y^d$ is a connected, $d$-dimensional *topological* manifold and $\phi : Y \to \partial_\infty X$ is a local homeomomorphism. The family of asymptotic Plateau problems forms a category where the morphisms between two asymptotic Plateau problems $(Y, \phi)$ and $(Y', \phi')$ are those continuous functions $\alpha : Y' \to Y$ such that

$$\phi' = \phi \circ \alpha. \tag{2.10}$$

In the case of $(d+1)$-dimensional hyperbolic space, asymptotic Plateau problems in $\partial_\infty \mathbb{H}^{d+1}$ are simply developed flat conformal structures, which have already been studied extensively in the literature (see, for example, our review [37]).

For $k \in \mathbb{N}$, we will say that an asymptotic Plateau problem is of **finite type** $k$ whenever it is isomorphic to a cover of $(\partial_\infty X \setminus F, \text{Id})$ for some finite subset $F$ of $\partial_\infty X$ of cardinality $k$. We obtain the following useful monotonicity result.

**Lemma 2.2.1**

*Let $X$ be a Cartan–Hadamard manifold with ideal boundary $\partial_\infty X$ and let $(Y, \phi)$ be an asymptotic Plateau problem in $\partial_\infty X$. If $(Y, \phi)$ contains an asymptotic Plateau problem of finite type $k$ for $k \leqslant 2$, then $(Y, \phi)$ is of finite type $l$, for some $l \leq k$.*

**Proof:** Let $U \subseteq Y$ be such that $(U, \phi)$ is of finite type $k$. If $k = 0$ then $U$ is compact and therefore also closed. Thus, since $Y$ is connected, $U = Y$, and the result follows trivially in this case.

Suppose now that $k > 0$. Let $F \subseteq \partial_\infty X$ be the complement of $\phi(U)$ in $\partial_\infty X$ which, by definition, has cardinality $k \leqslant 2$. Let $\partial U$ denote the topological boundary of $U$ in $Y$. We first show that

$$\phi(\partial U) \subseteq F. \tag{2.11}$$

Indeed, suppose the contrary. Let $y$ be a point of $\partial U$ whose image under $\phi$ is not an element of $F$, and let $V$ be a connected neighbourhood of $y$ in $Y$ over which $\phi$ restricts to a homeomorphism onto its image such that $\phi(V) \cap F = \emptyset$. Since the restriction of $\phi$ to $U$ is a covering of $\partial_\infty X \setminus F$, upon reducing $V$ if necessary, we may suppose that

$$\phi^{-1}(\phi(V)) \cap U = \sqcup_{\alpha \in A} W_\alpha,$$

where, for each $\alpha$, $\phi$ restricts to a homeomorphism from $W_\alpha$ onto $\phi(V)$. However, since $y \in \partial U$, $V \cap W_\alpha \neq \emptyset$ for some $\alpha$, and it follows by connectedness that $V = W_\alpha \subseteq U$, which is absurd, thus proving (2.11).

Now let $y$ be a point of $\partial U$ and let $V$ be a neighbourhood of $y$ in $Y$ homeomorphic to a ball over which $\phi$ restricts to a homeomorphism onto its image and $\phi(V) \cap F = \{\phi(y)\}$. We now show that $\phi$ restricts to a covering of $U \cup V$ onto its image. Indeed, observe first that $U \cap V$ is an open and closed subset of $V \setminus \{y\}$ so that, by connectedness, $V \setminus \{y\} \subseteq U$. It thus suffices to prove that $\phi$ is a covering near $\phi(y)$. However, since $k \leqslant 2$, $\phi$ maps the fundamental group of $V \setminus \{y\}$ surjectively onto the fundamental group of $\phi(U)$. In other words, every closed curve $\gamma$ in $\phi(U)$ starting and ending at some point $x$ of $\phi(V)$ is homotopic to $(\phi \circ \mu)$, for some curve $\mu$ contained in $V$, and the covering property follows. In particular, $U \cup V$ is of finite type $(k-1)$, and the result follows upon repeating this argument if necessary. $\square$

**2.3 - The Kulkarni–Pinkall section.** Let $X$ be a Cartan–Hadamard manifold with sectional curvature bounded above by $-1$. For any open horoball $H$ in $X$, we define its **asymptotic centre** $c(H)$ to be the unique point of intersection of its closure in $X \cup \partial_\infty X$ with $\partial_\infty X$. This function trivially makes the space $\text{HB}(X)$ of open horoballs in $X$ into a topological line bundle over $\partial_\infty X$ whose fibre over any point is the set of open horoballs asymptotically centred on that point. Furthermore, the relation of set inclusion defines a natural total order over every fibre of this bundle, where the empty set $\emptyset$ and the entire space $X$ identify with the end points of each fibre at $-\infty$ and $+\infty$ respectively.

Let $(Y, \phi)$ be an asymptotic Plateau problem in its ideal boundary $\partial_\infty X$. We define an **asymptotic ball** in $(Y, \phi)$ to be a pair $(\text{AB}, \alpha)$, where AB is an asymptotic ball in $\partial_\infty X$ and $\alpha : \text{AB} \to Y$ is a continuous open map such that

$$\phi \circ \alpha = \text{Id}. \tag{2.12}$$

Since $\phi$ is a local homeomorphism, every point of $Y$ is contained in some asymptotic ball.

Let $\phi^*\text{HB}(X)$ denote the pull-back of $\text{HB}(X)$ through $\phi$. In other words $\phi^*\text{HB}(X)$ is the topological line bundle over $Y$ whose fibre at any point $y$ is the set of open horoballs in $X$ asymptotically centred on $\phi(y)$. The section $\omega_\phi$ of this bundle is now constructed as follows. First, for $y \in Y$, we will say that an open horoball $H$ in $X$ asymptotically centred on $\phi(y)$ has **finite elevation** over $y$ whenever $y$ is contained in an asymptotic ball of the form $(\text{AB}(\xi), \alpha)$, for some inward-pointing unit normal vector $\xi$ of $\partial H$ at some point (see Figure 2.3.2). It is straightforward to see that, for all $y$, the set $\text{FE}_\phi(y)$ of open horoballs of finite elevation over $y$ is always non-empty. We then define

$$\omega_\phi(y) := \text{Sup}\, \text{FE}_\phi(y), \tag{2.13}$$

where, by convention,

$$\omega_\phi(y) := +\infty, \tag{2.14}$$

whenever this supremum does not exist. In the case where $X$ is $(d+1)$-dimensional hyperbolic space, we show in [37] that this section identifies with the area form of the metric constructed by Kulkarni–Pinkall in [20]. For this reason, we will call $\omega_\phi$ the **Kulkarni–Pinkall section** of $\phi$.

Kulkarni–Pinkall sections are monotone with respect to inclusion in the following sense.

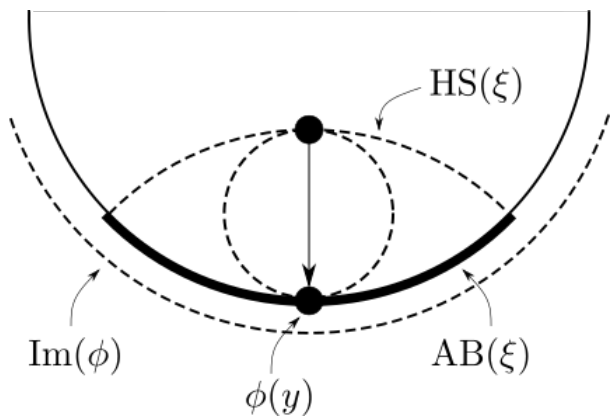

**Figure 2.3.2 - A horoball of finite elevation -** There is a point on the boundary of the horoball from which every tangent geodesic ray "lands" within the image of $\phi$.

**Lemma 2.3.1**

Let $X$ be a Cartan–Hadamard manifold with sectional curvature bounded above by $-1$ and let $\partial_\infty X$ denote its ideal boundary. Let $(Y, \phi)$ and $(Y', \phi')$ be asymptotic Plateau problems in $\partial_\infty X$ with respective Kulkarni–Pinkall sections $\omega_\phi$ and $\omega_{\phi'}$. For every morphism $\psi : Y' \to Y$,

$$\omega_{\phi'} \leqslant \psi^* \omega_\phi. \tag{2.15}$$

**Proof:** Indeed, if $y \in Y'$ lies in the asymptotic ball $(AB, \alpha)$, then trivially $\psi(y)$ lies in the asymptotic ball $(AB, \psi \circ \alpha)$. From this it follows that any horosphere in $X$ with finite elevation over $y$ also has finite elevation over $\psi(y)$, so that $FE_\phi(y) \subseteq FE_{\phi'}(\psi(y))$, and the result follows upon taking suprema. $\square$

The following finiteness result will play a key role in what follows.

**Lemma 2.3.2**

Let $X$ be a Cartan–Hadamard manifold with sectional curvature bounded above by $-1$ and let $\partial_\infty X$ denote its ideal boundary. Let $(Y, \phi)$ be an asymptotic Plateau problem in $\partial_\infty X$ and let $\omega_\phi$ denote its Kulkarni–Pinkall section.

(1) If $(Y, \phi)$ is of finite type 0 or 1, then $\omega_\phi$ is everywhere infinite.

(2) If $(Y, \phi)$ is not of finite type 0 or 1, then $\omega_\phi$ is everywhere finite.

**Remark 2.3.1.** The reader may verify that $\omega_\phi$ is in fact continuous, although we will have no need for this in the sequel.

**Proof:** The first assertion is trivial. To prove the second assertion, suppose the contrary. Thus, $(Y, \phi)$ is not of finite type 0 or 1 but there exists a point $y \in Y$ at which $\omega_\phi(y) = \infty$. Let $(H_m)_{m \in \mathbb{N}}$ be a sequence of open horoballs in $X$ asymptotically centred on $\phi(y)$ and exhausting $X$. For all $m$, let $(AB(\xi_m), \alpha_m)$ be an asymptotic ball in $(Y, \phi)$ containing $y$ where $\xi_m$ is an inward-pointing unit normal to $\partial H_m$ at some point. For all $m$, let $K_m$ denote the complement of $AB(\xi_m)$ in $\partial_\infty X$. Since the sectional curvature of $X$ is bounded above by $-1$, upon extracting a subsequence, we may suppose that $(K_m)_{m \in \mathbb{N}}$ converges in the Hausdorff sense to a single point $x$ of $\partial_\infty X$. In this case, the functions $(\alpha_m)_{m \in \mathbb{N}}$ join to define a local homeomorphism $\alpha : \partial_\infty X \setminus \{x\} \to Y$ such that $\phi \circ \alpha = \mathrm{Id}$. In particular, $(Y, \phi)$ contains an asymptotic Plateau problem of finite type 1, which is absurd by Lemma 2.2.1. This completes the proof. $\square$

**2.4 - Cartan–Hadamard ends.** Cartan–Hadamard ends play the same role in our theory as the hyperbolic ends constructed by Kulkarni–Pinkall in [20] do in the constant curvature case. They are defined as follows. First, we define a **height function** over a non-positively curved riemannian manifold $X$ to be a locally strictly convex $C^{1,1}_{\mathrm{loc}}$ function $h : X \to ]0, \infty[$ such that

(1) its gradient flow lines are unit speed geodesics; and

(2) for all $t > 0$, the superlevel set $h^{-1}([t, \infty[)$ is complete.

We will see in Corollary 2.4.4, below, that height functions, whenever they exist, are unique. We now define a **Cartan–Hadamard end** to be a non-positively curved riemannian manifold which carries a height function. This extends to the non-constant curvature case the concept of hyperbolic ends introduced by Kulkarni–Pinkall in [20] (see Section 3 of our review [37]). The family of all Cartan–Hadamard ends then forms a category whose morphisms are local isometries.

The simplest examples of Cartan–Hadamard ends are the complements of closed convex subsets of Cartan–Hadamard manifolds, where the height function is given by distance to the convex subset. In particular, for any Cartan–Hadamard manifold $X$ of sectional curvature pinched between $-K$ and $-1$ for some $K \geqslant 1$, and for any unit vector $\xi$ in its unit sphere bundle, the canonical end $CE(\xi)$, constructed in Section 2.1, is a Cartan–Hadamard end with height function given by distance to the boundary. This observation will prove useful later.

Cartan–Hadamard ends possess considerable geometric structure which it is worth describing in some detail. Let $X$ be a Cartan–Hadamard end with height function $h$ and let $TX$ denote its tangent bundle. We say that a tangent vector $\xi \in TX$ over some point $x$ is **upward pointing** whenever

$$\langle \xi, \nabla h(x) \rangle > 0. \tag{2.16}$$

We denote respectively by $T^+X$ and $S^+X$ the subbundles of upward pointing tangent vectors and unit upward pointing tangent vectors over $X$. We call the gradient flow lines of $h$ **vertical lines**. These curves form a geodesic foliation of $X$ which we call its **vertical line foliation** and whose leaf space we denote by $\mathcal{V}$. We call the level sets of $h$ the **levels** of $X$. These form another foliation of $X$ by $C^{1,1}_{\text{loc}}$ embedded hypersurfaces which we call the **level set foliation** and which we denote by $(X_t)_{t>0}$. These two foliations are transverse to one another and every vertical line intersects every level at exactly one point, from which it follows that every level is naturally homeomorphic to $\mathcal{V}$. For all $t > 0$, the vertical projection $\pi_t : X \to X_t$ is defined to be the function which sends each point $x$ of $X$ to the unique intersection with $X_t$ of the vertical line in which it lies. By standard properties of convex subsets of non-positively curved manifolds, for all $t > 0$, $\pi_t$ restricts to a 1-Lipschitz function from $h^{-1}([t, \infty[)$ onto $X_t$.

The following analogue of the Hopf-Rinow theorem encapsulates the main geometric properties of Cartan–Hadamard ends.

### Theorem 2.4.1, Hopf-Rinow

*Let $X$ be a Cartan–Hadamard end with height function $h$. For all $x \in X$, every tangent vector $\xi \in T_xX$ satisfying*

$$\langle \xi_x, \nabla h(x) \rangle > -h(x) \tag{2.17}$$

*lies in the domain of $\mathrm{Exp}_x$. In particular $T^+X$ is contained in the domain of $\mathrm{Exp}$.*

**Proof:** Indeed, by convexity, for all $t \in [0, 1]$ such that $t\xi_x$ lies in the domain of $\mathrm{Exp}_x$,

$$(h \circ \mathrm{Exp}_x)(t\xi_x) \geqslant h(x) + \langle \xi_x, \nabla h(x) \rangle t,$$

so that,

$$\mathrm{Exp}_x(t\xi_x) \in h^{-1}([h(x) + \langle \xi_x, \nabla h(x) \rangle, \infty[),$$

and the result now follows by completeness. $\square$

### Corollary 2.4.2

*Let $X$ be a Cartan–Hadamard end with height function $h$. If $\gamma : [0, \alpha[ \to X$ is a unit speed geodesic segment such that*

$$\langle \dot{\gamma}(0), (\nabla h \circ \gamma)(0) \rangle \geqslant 0. \tag{2.18}$$

*then $\gamma$ extends to a geodesic ray defined over the entire half-line $[0, \infty[$.*

**Corollary 2.4.3**

*Let $X$ be a Cartan–Hadamard end with height function $h$. For all $x \in X$, the only point of the closed ball of radius $h(x)$ about $0$ in $T_xX$ which does not lie in the domain of $\mathrm{Exp}_x$ is $-h(x)\nabla h(x)$.*

**Corollary 2.4.4**

*Let $X$ be a Cartan–Hadamard end with height function $h$. For all $x \in X$,*

$$h(x) := \mathrm{Sup}\left\{r \mid \overline{B}_r(x) \text{ is compact}\right\}. \tag{2.19}$$

*In particular, the height function, when it exists, is unique.*

We conclude this section by studying the geometric structures of morphisms between Cartan–Hadamard ends. We first show that the height function is monotone with respect to inclusion.

**Lemma 2.4.5**

*Let $X$ and $X'$ be Cartan–Hadamard ends with respective height functions $h$ and $h'$. If $\phi : X' \to X$ is a local isometry, then*

$$h \circ \phi \geqslant h'. \tag{2.20}$$

**Proof:** Indeed, choose $x \in X'$ and denote $y := \phi(x)$. Since $\phi$ is a local isometry,

$$\mathrm{Exp}_y \circ D\phi(x) = \phi \circ \mathrm{Exp}_x.$$

Consequently, for all $r > 0$, if $\overline{B}_r(0)$ is contained in the domain of $\mathrm{Exp}_x$, then $D\phi(x)(\overline{B}_r(0))$ is contained in the domain of $\mathrm{Exp}_y$. The result now follows by (2.19). $\square$

Finally, we will say that an embedded hypersurface in a Cartan–Hadamard end is a **graph** whenever it is transverse to the vertical line foliation and meets each vertical line at most once.

**Lemma 2.4.6**

*Let $X$ and $X'$ be Cartan–Hadamard ends and let $\phi : X' \to X$ be an injective local isometry. If $(X'_t)_{t>0}$ denotes the level set foliation of $X'$ then, for all $t$, $\phi(X'_t)$ is a graph in $X$.*

**Proof:** Indeed, let $\gamma : ]0,\infty[ \to X$ be a height-parametrised vertical line. Suppose that $\gamma$ meets $\alpha(X'_t)$ at some point, and let $t_0 > 0$ be the first height at which this occurs. Denote $x := \phi(y) := \gamma(t_0)$. Let $\mu : ]a,\infty[ \to X'$ be the unique geodesic such that
$$D\phi(y)\dot\mu(0) = \dot\gamma(t_0).$$

We now show that $\gamma$ is transverse to $\phi(X'_t)$ at $t_0$. Indeed, otherwise, $\mu$ is tangent to $X'_t$ at $y$, so that, by Corollary 2.4.2, $\mu$ extends to a geodesic defined over the whole of $\mathbb{R}$. Since $\gamma = \phi \circ \mu$, it follows that $\gamma$ also extends to a geodesic defined over the whole of $\mathbb{R}$. This is absurd, since $\gamma$ is parametrised by height, and transversality follows. It now remains only to show that $\gamma$ does not intersect $\phi(X'_t)$ at any other point. However, let $h'$ denote the height function of $X'$. By convexity, $(h' \circ \mu)$ is strictly increasing over $]0,\infty[$. Consequently,

$$\mu(]0,\infty[) \cap X'_t = \emptyset,$$

so that, by injectivity,

$$\gamma(]t_0,\infty[) \cap \phi(X'_t) = \phi(\mu(]0,\infty[) \cap X'_t) = \emptyset.$$

The only point of intersection of $\gamma$ with $\phi(X'_t)$ is thus indeed $\gamma(t_0)$, as asserted, and this completes the proof. $\square$

**2.5 - Local analytic properties of the height function.** We now study in more detail the analytic structures of height functions. For this we recall the following lesser used definition from the theory of partial differential equations. Given a smooth manifold $M$, a point $p \in M$, a function $f : M \to \mathbb{R}$ and a symmetric bilinear form $B$ over $T_pM$, we say that

$$\text{Hess}(f)(p) \geq B, \tag{2.21}$$

in the **weak sense** whenever, for all $\epsilon > 0$, there exists a neighbourhood $\Omega$ of $p$ in $M$ and a smooth function $g : \Omega \to \mathbb{R}$ such that

(1) $g \leqslant f$;

(2) $g(p) = f(p)$; and

(3) $\text{Hess}(g)(p) = B - \epsilon \text{Id}$.

We likewise say that $\text{Hess}(f)(p) \leqslant B$ in the **weak sense** whenever $\text{Hess}(-f)(p) \geqslant -B$ in the weak sense. This concept was first used by Calabi in [8] and the reader will notice a similarity with the more modern concept of viscosity solutions (see [12]). However, we must emphasize that Calabi's concept is more restrictive, since the theory of viscosity solutions does not actually require that the comparison functions *exist* at every point.

We now obtain weak estimates of the hessian of the height function. We refer the reader to Appendix A for a brief review of the terminology used in the statement of the following theorem.

**Lemma 2.5.1**

*Let $X$ be a Cartan–Hadamard end with height function $h$.*

*(1) If the sectional curvature of $X$ is bounded above by $-k \leq 0$ then, for all $x \in X$, with respect to the decomposition $T_xX = \text{Ker}(dh(x)) \oplus \langle \nabla h(x) \rangle$,*

$$\text{Hess}(h)(x) \geqslant \begin{pmatrix} \tanh_k(h(x))\text{Id} & 0 \\ 0 & 0 \end{pmatrix}, \tag{2.22}$$

*in the weak sense.*

*(2) If the sectional curvature of $X$ is bounded below by $-K \leq 0$ then, for all $x \in X$, with respect to the decomposition $T_xX = \text{Ker}(dh(x)) \oplus \langle \nabla h(x) \rangle$,*

$$\text{Hess}(h)(x) \leqslant \begin{pmatrix} \coth_K(h(x))\text{Id} & 0 \\ 0 & 0 \end{pmatrix}, \tag{2.23}$$

*in the weak sense.*

**Remark 2.5.1.** In particular, when the sectional curvature of $X$ is bounded below, $\text{Hess}(h)$ is uniformly bounded over $h^{-1}([t,\infty[)$ in the weak sense for all $t > 0$.

Before proving 2.5.1, we observe the following useful corollary of this result and (B.1).

**Corollary 2.5.2**

*Let $X$ be a Cartan–Hadamard end with sectional curvature pinched between $-K$ and $-1$. For all $t > 0$, the first fundamental form $I_t$ and the second fundamental form $II_t$ of the level $X_t$ satisfy*

$$\tanh(t)I_t \leqslant II_t \leqslant \tanh_K(t)I_t \tag{2.24}$$

*in the weak sense.*

**Proof of Lemma 2.5.1:** Choose $x \in X$ and $t < h(x)$. Let $y \in X$ be the point lying vertically below $x$ at height $h(x) - t$. Define the subset $E \subseteq T_yX$ by

$$E := \{\xi_y \mid \langle \xi_y, \nabla h(y) \rangle = 0\}.$$

By Theorem 2.4.1, $E$ is contained in the domain of the exponential map of $X$. Let $d_y$ and $d_E$ denote the respective distances in $X$ to the point $y$ and the hypersurface $\text{Exp}_y(E)$. Trivially, $d_E(x) = h(x) - h(y)$. Furthermore, by convexity, $h$ achieves its minimum over $E$ at $y$ so that, by Corollary 2.4.4, $d_E \leqslant h - h(y)$. Likewise, $d_y(x) = h(x) - h(y)$ and, by Corollaries 2.4.3 and 2.4.4, $d_y \geqslant h - h(y)$. It follows that, at $x$,

$$\text{Hess}(d_E)(x) \leqslant \text{Hess}(h)(x) \leqslant \text{Hess}(d_y)(x),$$

in the weak sense, and the result now follows by classical comparison theory upon letting $r$ tend to $h(x)$. $\square$

**2.6 - The ideal boundary of a Cartan–Hadamard end.** We use the estimates of the preceding sections to study the ideal boundaries of Cartan–Hadamard ends. We first divide complete, unit speed geodesic rays in Cartan–Hadamard ends into two classes.

**Lemma and Definition 2.6.1**

*Let $X$ be an Cartan–Hadamard end with sectional curvature bounded above by $-1$. For every complete, unit speed geodesic ray $\gamma : [0, \infty[ \to X$,*

$$\operatorname*{Lim}_{t \to \infty}(h \circ \gamma)(t) \in \{0, \infty\}. \tag{2.25}$$

*We say that $\gamma$ is **bounded** whenever this limit is zero, and **unbounded** otherwise.*

**Proof:** Indeed, by convexity $(h \circ \gamma)(t)$ converges to a (possibly infinite) limit as $t$ tends to infinity. Suppose now that
$$\lim_{t \to +\infty}(h \circ \gamma)(t) > 2\epsilon,$$
for some $\epsilon > 0$. Denoting $f := (h \circ \gamma)(t)$, (2.22) yields, for sufficiently large $t$,

$$\ddot{f} \geqslant \tanh(\epsilon)(1 - \dot{f}^2), \tag{2.26}$$

in the weak sense. Upon solving this ordinary differential inequality, we see that $f(t)$ tends to $+\infty$ as $t$ tends to infinity, as desired. $\square$

**Lemma 2.6.2**

*Let $X$ be a Cartan–Hadamard end with sectional curvature bounded above by $-1$, and let $\gamma : [0, \infty[ \to X$ be an unbounded geodesic ray in $X$. For all $t > 0$, there exists $x \in X_t$ such that*

$$\operatorname*{Lim}_{s \to +\infty}(\pi_t \circ \gamma)(s) = x. \tag{2.27}$$

*In particular, $\gamma$ is asymptotic to the vertical line passing through $x$.*

**Proof:** Indeed, let $h$ denote the height function of $X$ and denote $f := (h \circ \gamma)$. Since $\gamma$ is unbounded, by (2.22), for all $\epsilon > 0$ and for all sufficiently large $t$,

$$\ddot{f} \geqslant (1 - \epsilon)(1 - \dot{f}^2),$$

in the weak sense. Upon solving this ordinary differential inequality, we see that there exists $a > 0$ such that, for all sufficiently large $t$,

$$\langle \dot{\gamma}(t), (\nabla h \circ \gamma)(t) \rangle = \dot{f}(t) \geqslant \tanh((1-\epsilon)(t-a)).$$

In particular, the orthogonal projection of $\dot{\gamma}(t)$ onto $\langle \nabla h \circ \gamma)(t) \rangle^\perp$ decays exponentially as $t$ tends to infinity. Since $\pi_t$ is 1-Lipschitz, and since $\operatorname{Ker}(D\pi_t(x)) = \langle \nabla h(x) \rangle$ for all $x$, it follows that the projection $(\pi_t \circ \gamma)$ has finite length in $X_t$. The result now follows by completeness. $\square$

Let $X$ be a Cartan–Hadamard end of sectional curvature bounded above by $-1$. We define $\partial_\infty X$, the **ideal boundary** of $X$, to be the space of equivalence classes of unbounded, unit speed geodesic rays in $X$, where two such rays are deemed equivalent whenever they are asymptotic to one another. We furnish $X \cup \partial_\infty X$ with the cone topology, defined in an analogous manner to that of the ideal boundary of a Cartan–Hadamard manifold, as described in Section 2.1.

By Lemma 2.6.2, every unbounded geodesic ray in $X$ is asymptotic to some vertical line. On the other hand, since $\pi_t$ is 1-Lipschitz for all $t$, no two vertical lines are asymptotic to one another. It follows that $\partial_\infty X$ is homeomorphic to the leaf space $\mathcal{V}$ of the vertical line foliation which, we recall, is in turn homeomorphic to every level $X_t$ of $X$. We define $\pi_\infty : X \to \partial_\infty X$ to be the function which sends every point of $x$ to the equivalence class of the vertical line on which it lies. By the preceding discussion, for all $t$, $\pi_\infty$ restricts to a homeomorphism from $X_t$ into $\partial_\infty X$.

On the asymptotic Plateau problem in Cartan-Hadamard manifolds.

Recall that $S_+X$ denotes the bundle of upward pointing unit tangent vectors over $X$. Bearing in mind Corollary 2.4.2, as in Section 2.1, for all $\xi \in S_+X$, we define the geodesic ray $\gamma_\xi : [0, \infty[ \to X$ such that

$$\dot{\gamma}_\xi(0) = \xi. \qquad (2.28)$$

We then define the **horizon map** $\text{Hor} : S_+X \to \partial_\infty X$ by

$$\text{Hor}(\xi) := \lim_{t \to \infty} \gamma_\xi(t). \qquad (2.29)$$

We verify as before that this function is continuous and restricts to a homeomorphism of every fibre of $S_+X$ onto its image.

We conclude this section by reviewing the functorial properties of the ideal boundary operator. Thus, let $X$ and $X'$ be Cartan–Hadamard ends with sectional curvatures bounded above by $-1$ and let $\phi : X' \to X$ be a local isometry. Let $\gamma : [0, \infty[ \to X'$ be a complete, unit-speed geodesic ray. By Lemma 2.4.5, if $\gamma$ is unbounded, then so too is $\phi \circ \gamma$. Thus, since $\phi$ sends equivalent geodesic rays to equivalent geodesic rays, it induces a function $\partial_\infty \phi : \partial_\infty X' \to \partial_\infty X$ which, by standard properties of convex subsets of Cartan–Hadamard manifolds, is a local homeomorphism. In addition, this operator respects the identity in the sense that, if $\phi := \text{Id}$, then $\partial_\infty \phi = \text{Id}$, and it respects compositions in the sense that, if $X$, $X'$ and $X''$ are Cartan–Hadamard ends, and if $\phi' : X'' \to X'$ and $\phi : X' \to X$ are local isometries, then

$$\partial_\infty(\phi \circ \phi') = (\partial_\infty \phi) \circ (\partial_\infty \phi'). \qquad (2.30)$$

In summary, $\partial_\infty$ defines a covariant functor from the category of Cartan–Hadamard ends into the category of topological manifolds.

**2.7 - Developed Cartan–Hadamard ends.** Let $X$ be a Cartan–Hadamard manifold with sectional curvature bounded above by $-1$ and let $\partial_\infty X$ denote its ideal boundary. We define a **developed Cartan–Hadamard end** in $X$ to be a pair $(Y, \phi)$ where $Y$ is a Cartan–Hadamard end and $\phi : Y \to X$ is a local isometry. The family of developed Cartan–Hadamard ends in $X$ also forms a category, where the morphisms of interest to us between two developed Cartan–Hadamard ends $(Y, \phi)$ and $(Y', \phi')$ are those injective local isometries $\alpha : Y \to Y'$ such that $\phi = \phi' \circ \alpha$. We verify, as in the preceding section, that $\phi$ induces a local homeomorphism $\partial_\infty \phi : \partial_\infty Y \to \partial_\infty X$ which respects compositions. That is, $(\partial_\infty Y, \partial_\infty \phi)$ is an asymptotic Plateau problem in $\partial_\infty X$. It follows that $\partial_\infty$ defines a covariant functor from the category of developed Cartan–Hadamard ends in $X$ into the category of asymptotic Plateau problems in its ideal boundary.

Let $(Y, \phi)$ be a developed Cartan–Hadamard end in $X$ with height function $h$. We define an **open half-space** in $Y$ to be a pair $(\text{HS}, \alpha)$ where HS is an open half-space in $X$ and $\alpha : \text{HS} \to Y$ is a local isometry such that

$$\phi \circ \alpha = \text{Id}. \qquad (2.31)$$

We construct a large family of open half spaces in $Y$ as follows. Recall that $T_+Y$ denotes the bundle of upward pointing tangent vectors over $Y$. Choose $y \in Y$ and denote $x := \phi(y)$ and $\xi := D\phi(y) \cdot \nabla h(y)$. Observe that $\text{Exp}_x \circ D\phi(y)$ maps $T_{+,y}Y$ diffeomorphically onto the open half-space $\text{HS}(\xi)$. We thus denote

$$\text{HS}(y) := \text{HS}(\xi), \qquad (2.32)$$

and we define $\alpha(y) : \text{HS}(y) \to Y$ such that

$$\alpha(y) \circ \text{Exp}_x \circ D\phi(y) = \text{Exp}_y. \qquad (2.33)$$

We call $(\text{HS}(y), \alpha(y))$ the **open half-space** in $Y$ **centred** on $y$. We likewise denote

$$\text{CE}(y) := \text{CE}(\xi), \qquad (2.34)$$

and, for all $t > 0$, we denote

$$\text{HS}_t(y) := \text{HS}_t(\xi). \qquad (2.35)$$

Recall that, for all $y$, $CE(y)$ is itself a Cartan–Hadamard end with height function given by distance to the boundary. In particular, for all $y$, $\alpha(y)$ defines a morphism from this Cartan–Hadamard end into $Y$.

We likewise construct a large family of asymptotic balls in $\partial_\infty Y$ as follows. Let $y$, $x$ and $\xi$ be as before. Let $\mathrm{Hor}^X$ and $\mathrm{Hor}^Y$ denote the respective horizon maps of $X$ and $Y$. Observe that $\mathrm{Hor}^X \circ D\phi(y)$ maps $S_{+,y}Y$ homeomorphically onto $\mathrm{AB}(\xi)$. We thus denote

$$\mathrm{AB}(y) := \mathrm{AB}(\xi), \tag{2.36}$$

and we define $\alpha(y) : \mathrm{AB}(y) \to \partial_\infty Y$ such that

$$\alpha(y) \circ \mathrm{Hor}^X_x \circ D\phi(y) = \mathrm{Hor}^Y_y. \tag{2.37}$$

We call $(\mathrm{AB}(y), \alpha(y))$ the **asymptotic ball** in $\partial_\infty Y$ **centred** on $y$. Observe, in particular, that this ball contains $\pi_\infty(y)$.

The following estimate will prove useful in the study of the non-complete case.

**Lemma 2.7.1**

Let $X$ be a Cartan–Hadamard manifold with sectional curvature bounded above by $-1$. Let $(Y, \phi)$ be a developed Cartan–Hadamard end in $X$ with level set foliation $(Y_t)_{t>0}$ and ideal boundary $\partial_\infty Y$. For all $t$, let $\pi_t : \partial_\infty Y \to Y_t$ denote the canonical projection. For all $y \in Y$ and for all $0 < t < h(y)$,

$$\mathrm{Diam}((\pi_t \circ \alpha(y))(\mathrm{AB}(y))) \leqslant 2\mathrm{arccosh}(\mathrm{coth}(h(y) - t)). \tag{2.38}$$

**Proof:** By convexity and comparison theory, it is sufficient to consider the model case where $X$ is 2-dimensional hyperbolic space, $Y$ is an open half-space in $X$ and $\pi := \pi_0$ is the orthogonal projection onto the geodesic $\Gamma := \partial Y$. The result is then an exercise of classical hyperbolic geometry. $\square$

**2.8 - Infinitesimally strictly convex immersed hypersurfaces.** We conclude this chapter by associating developed Cartan–Hadamard ends in a natural manner to quasicomplete ISC immersed hypersurfaces in Cartan–Hadamard manifolds. In particular, this will yield a new $C^0$ a priori estimate for such hypersurfaces which will play a fundamental role in our proof of existence of solutions to the asymptotic Plateau problem. Let $X$ be a Cartan–Hadamard manifold and let $(Y, e)$ be an ISC immersed hypersurface in $X$. Let $\nu_e : Y \to SX$ denote the positively oriented unit normal vector field over $Y$, and define

$$\begin{aligned} \hat{Y} &:= Y \times ]0, \infty[, \text{ and} \\ \hat{e}(x,t) &:= \mathrm{Exp}(t\nu_e(x)). \end{aligned} \tag{2.39}$$

By standard properties of ISC immersions in Cartan–Hadamard manifolds, $\hat{e}$ is a local diffeomorphism, and we furnish $\hat{Y}$ with the unique riemannian metric that makes $\hat{e}$ into a local isometry.

**Lemma and Definition 2.8.1**

Let $X$ be a Cartan–Hadamard manifold with sectional curvature pinched between $-K$ and $-1$, and let $(Y, e)$ be an ISC immersed hypersurface in $X$. $(\hat{Y}, \hat{e})$ is a developed Cartan–Hadamard end in $X$ if and only if $(Y, e)$ is quasicomplete. When this holds, we call $(\hat{Y}, \hat{e})$ the **end** of $(Y, e)$ and its height function is given by

$$h(x, t) := t. \tag{2.40}$$

**Proof:** By standard properties of ISC immersions in Cartan–Hadamard manifolds, $h(x,t) := t$ is a convex $C^{1,1}_{\mathrm{loc}}$ function whose gradient flow lines are unit speed geodesics. It thus remains only to study completeness. However, for all $t$, denote $e_t(x) := \hat{e}(x,t)$ and let $\mathrm{I}_t$, $\mathrm{II}_t$ and $\mathrm{III}_t$ denote respectively the first, second and third fundamental forms of this immersion. By classical comparison theory, for all $t$,

$$\cosh^2(t)\mathrm{I}_e + 2\cosh(t)\sinh(t)\mathrm{II}_e + \sinh^2(t)\mathrm{III}_e \leqslant \mathrm{I}_t \leqslant \cosh^2_K(t)\mathrm{I}_e + 2\cosh_K(t)\sinh_K(t)\mathrm{II}_e + \sinh^2_K(t)\mathrm{III}_e,$$

where the functions $\cosh_K$ and $\sinh_K$ are defined in Appendix A. Since $\mathrm{II}_e \geqslant 0$, this yields

$$\sinh^2(t)(\mathrm{I}_e + \mathrm{III}_e) \leqslant \mathrm{I}_t \leqslant 2\cosh^2_K(t)(\mathrm{I}_e + \mathrm{III}_e). \tag{2.41}$$

It follows that $(Y, e)$ is quasicomplete if and only if $\mathrm{I}_t$ is complete for all $t$. However, since $\pi_t$ is 1-Lipschitz for all $t$, the latter holds if and only if $h^{-1}([t, \infty[)$ is complete for all $t$, and the result follows. $\square$

This yields the desired $C^0$ a priori estimate.

On the asymptotic Plateau problem in Cartan-Hadamard manifolds.**Theorem 2.8.2**

*Let $X$ be a Cartan–Hadamard manifold with sectional curvature pinched between $-K$ and $-1$ and let $\partial_\infty X$ denote its ideal boundary. Let $(Y,e)$ be a quasicomplete ISC immersed submanifold in $X$, let $\phi_e : Y \to \partial_\infty X$ denote its asymptotic Gauss map and let $\omega_e$ denote the Kulkarni–Pinkall section of $\phi_e$. For all $y \in Y$,*

$$e(y) \in \overline{\omega_e}(y). \tag{2.42}$$

**Remark 2.8.1.** In fact, it is straightforward to show that $e(y)$ is contained in the interior of this horoball, although (2.42) is already sufficient for our purposes.

**Proof:** Let $(\hat{Y}, \hat{e})$ denote the end of $(Y, e)$. Observe that $Y$ identifies with $\partial \hat{Y}$ and that $\pi_\infty : Y \to \partial_\infty \hat{Y}$ is a homeomorphism. Choose $(y,t) \in \hat{Y}$ and denote $x := \hat{e}(y,t)$ and $\xi := D\hat{e}(y,t) \cdot \partial_t$. Let $(\text{AB}(y,t), \alpha(y,t))$ denote the asymptotic ball of $\partial_\infty \hat{Y}$ centred on $(y,t)$. Let $H(y,t)$ denote the unique open horosphere in $X$ with boundary passing through $x$ and inward-pointing normal $\xi$ at this point. Since

$$\text{AB}(y,t) = \text{AB}(\xi),$$

$H(y,t)$ has finite elevation over $\pi_\infty(y,t)$. Consequently,

$$\hat{e}(y,t) = x \in \partial H \subseteq \overline{\omega}_{\partial_\infty \hat{e}}(\pi_\infty(y,t)).$$

Upon letting $t$ tend to zero, this yields,

$$e(y) = \lim_{t \to 0} \hat{e}(y,t) \in \overline{\omega}_{\partial_\infty \hat{e}}(\pi_\infty(y)).$$

Finally, for all $z \in Y$,

$$(\partial_\infty \hat{e} \circ \pi_\infty)(z) = \lim_{s \to \infty} \hat{e}(z,s) = \lim_{s \to \infty} \text{Exp}(t\nu_e(z)) = \phi_e(z),$$

so that

$$\omega_e(y) = \omega_{\partial_\infty \hat{e}}(\pi_\infty(y)),$$

and the result follows. $\square$

## 3 - Special lagrangian curvature, level sets and uniqueness.

**3.1 - Special lagrangian curvature.** In [35], we introduced the concept of special lagrangian curvature of ISC immersed hypersurfaces in riemannian manifolds. This curvature notion, which arises from the contact structure of the unit sphere bundle, extends to higher dimensions the remarkable properties of 2-dimensional extrinsic curvature studied by Labourie in a series of papers (see, for example, [21], [22] and [23]) and described in full generality in [24].

In order to construct special lagrangian curvature, we first define the special lagrangian angle of a real, symmetric matrix. First observe that, for every real, symmetric matrix $A$,

$$\text{Det}(\text{Id} + iA) \neq 0. \tag{3.1}$$

Consequently, since the space of real, symmetric matrices is contractible, there exists a unique, real-analytic function $\Theta : \text{Symm}(\mathbb{R}^d) \to \mathbb{R}$, which we call the **special lagrangian angle**, such that

$$\Theta(0) = 0, \tag{3.2}$$

and, for all $A$,

$$\text{Im}(\text{Log}(\text{Det}(I + iA))) = \Theta(A) \mod 2\pi. \tag{3.3}$$

It is straightforward to show that, for any real symmetric $A$ with eigenvalues $\lambda_1, \cdots, \lambda_d$,

$$\Theta(A) = \sum_{i=1}^{d} \arctan(\lambda_i), \tag{3.4}$$

so that, in particular, for all $A$,

$$\Theta(A) \in ]-d\pi/2, d\pi/2[. \tag{3.5}$$

**Lemma 3.1.1**

*The first and second derivatives of the special lagrangian angle are given by*

$$D\Theta(A) \cdot B = \mathrm{Tr}((Id + A^2)^{-1}B), \text{ and}$$
$$D^2\Theta(A)(B,B) = -2\mathrm{Tr}((Id + A^2)^{-1}AB(Id + A^2)^{-1}B). \tag{3.6}$$

*In particular, the special lagrangian angle is strictly concave over the cone of positive definite matrices.*

**Proof:** Indeed,

$$\mathrm{D}\Theta(A) \cdot B = \mathrm{Im}(i\mathrm{Tr}((I + iA)^{-1}B)) = \mathrm{Re}(\mathrm{Tr}((I - iA)(I + A^2)^{-1}B)) = \mathrm{Tr}((I + A^2)^{-1}B),$$

and the first relation follows. The second relation is trivially obtained upon differentiating a second time, and this completes the proof. □

Special lagrangian curvature is now defined in terms of the special lagrangian angle. However, it will first be useful to recall the properties that should be expected of curvature functions as well as their motivations (c.f. [36] and the references therein). To this end, let $\Lambda^d \subseteq \mathrm{End}(\mathbb{R}^d)$ denote the open cone of positive definite symmetric matrices and let $\overline{\Lambda}^d$ denote its closure. We define a **convex curvature function** to be a continuous function $K : \overline{\Lambda}^d \to \mathbb{R}$, differentiable over the interior of this cone, such that

(1) $K(M^t A M) = K(A)$ for all $A$ and for every orthogonal matrix $M$;

(2) $K$ vanishes over $\partial\Lambda^d$;

(3) $K$ is 1-homogeneous;

(4) $K(\mathrm{Id}) = 1$;

(5) $DK(A) \cdot B > 0$ for all $A \in \Lambda^d$ and for all non-negative symmetric $B$; and

(6) $K$ is concave.

The first four of these conditions ensure that $K$ defines a geometrically well-behaved scalar notion of curvature for ISC immersed hypersurfaces. Indeed, let $A_e$ denote the Weingarten operator of some ISC immersed hypersurface $(Y, e)$, say, in some riemannian manifold $X$. $O(d)$ invariance ensures that $K_e := K(A_e)$ is a well-defined function over $Y$ which we call its $K$-**curvature**. The vanishing of $K$ over $\partial\Lambda^d$ ensures that $K$-curvature is compatible with convexity, 1-homogeneity ensures that it transforms correctly under rescalings of the ambient metric, and the normalisation of Condition (4) ensures that the $K$-curvature of the unit sphere in euclidean space is equal to 1. The final two conditions are analytic conditions required for $K$-curvature to be studied via existing partial differential equations techniques. Condition (5), known as **ellipticity**, ensures that the Jacobi operator of $K$-curvature of any ISC immersed hypersurface is elliptic, whilst concavity is a standard condition used in the study of totally non-linear partial differential equations (c.f., for example, [7]).

We now proceed to the construction of a 1-parameter family of smooth convex curvature functions. First, for $r > 0$, denote

$$\Theta_r(A) := \Theta(A/r). \tag{3.7}$$

Observe that, for all $A \in \Lambda^d$, $\Theta_r$ is strictly decreasing in $r$ and

$$\lim_{r \to 0} \Theta_r(A) = d\pi/2, \text{ and}$$
$$\lim_{r \to \infty} \Theta_r(A) = 0. \tag{3.8}$$

Consequently, for all $\theta \in ]0, d\pi/2[$, there exists a unique function $R_\theta : \Lambda^d \to \mathbb{R}$ such that, for all $A \in \Lambda^d$,

$$\Theta_{R_\theta(A)/\tan(\theta/d)}(A) = \theta. \tag{3.9}$$

**Lemma 3.1.2**

*For all $\theta \in [(d-1)\pi/2, d\pi/2[$, $R_\theta$ is a convex curvature function which is smooth over $\Lambda^d$.*

**Proof:** Smoothness of $R_\theta$ follows by (3.6) and the implicit function theorem. By $O(d)$-invariance of the determinant, $\Theta$ and $R_\theta$ are also $O(d)$-invariant. For all $r > 0$ and for all $A \in \partial\Lambda^d$, $\Theta_r(A) < (d-1)\pi/2$, from which it follows that $R_\theta(A) = 0$ for all $\theta \geqslant (d-1)\pi/2$. By construction, $R_\theta$ is homogeneous of order 1 and $R_\theta(\mathrm{Id}) = 1$. By (3.6) and the implicit function theorem again, for all $A \in \Lambda^d$ and for all non-negative, symmetric $B$, $DR_\theta(A) \cdot B > 0$. It remains only to prove concavity of $R_\theta$. Since this function is homogeneous of order 1, it suffices to show that its superlevel sets are strictly convex. However, since $R_\theta(A) > r$ if and only if $\Theta_r(A) > \theta$, this property follows by concavity of $\Theta$, and this completes the proof. $\square$

In this paper, we will only be concerned with the case where $\theta = (d-1)\pi/2$, and we thus denote

$$R := R_{(d-1)\pi/2}, \tag{3.10}$$

We call $R$ the **special lagrangian curvature function**. Let $X$ be a riemannian manifold. Let $(Y, e)$ be an ISC immersed hypersurface in $X$, let $\nu_e$ denote its outward pointing unit normal vector field and let $A_e$ denote its corresponding Weingarten operator. We define the **special lagrangian curvature** of $e$ by

$$R_e := R(A_e), \tag{3.11}$$

and, for all $r > 0$, we define its $r$-**special Lagrangian angle** by

$$\Theta_{r,e} := \Theta_r(A_e). \tag{3.12}$$

Despite having a rather involved construction, in low dimensions special lagrangian curvature enjoys relatively simple expressions in terms of classical notions of curvature. Indeed, in the case where $X$ is 3-dimensional and $(Y, e)$ is a surface, the special lagrangian curvature is none other than the square root of the extrinsic curvature, that is

$$R_e = \sqrt{K_e}. \tag{3.13}$$

Likewise, in the case where $X$ is 4-dimensional, we still have the following relatively simple formula for the special lagrangian curvature

$$R_e = \sqrt{3\mathrm{Det}(A_e)/\mathrm{Tr}(A_e)}. \tag{3.14}$$

Nevertheless, the formula for $R_e$ in terms of classical curvature notions becomes progressively less trivial as the dimension grows. However, the main appeal of special lagrangian curvature is that it constitutes the correct higher-dimensional framework within which the theory developed by Labourie for constant extrinsic curvature surfaces continues to apply, as we will show in Section 4.1, below.

**3.2 - Foliations of Cartan–Hadamard ends.** We now study the special lagrangian angles of the levels of Cartan–Hadamard ends. This will allow us to prove in the following section a uniqueness result for solutions of asymptotic Plateau probems.

**Theorem 3.2.1**

*Let $X$ be a Cartan–Hadamard end with smooth boundary and with sectional curvature bounded above by $-1$ and let $(X_t)_{t>0}$ denote its level set foliation. For $(d-1)\pi/2 \leqslant \theta_0 < d\pi/2$ and $0 < r\tanh(\theta_0/d) \leqslant 1$, if $\partial X$ has $r$-special lagrangian angle no less than $\theta_0$ then, for all $t > 0$, $X_t$ has $r$-special lagrangian angle no less than*

$$\theta(t; t, \theta_0) := d\arctan(\tanh(t_0 + t)/r), \tag{3.15}$$

*where*

$$t_0 := \mathrm{arctanh}(r\tan(\theta_0/d)). \tag{3.16}$$

Theorem 3.2.1 is a consequence of the following estimate of [34]. We include the proof for the reader's convenience.

On the asymptotic Plateau problem in Cartan-Hadamard manifolds.

**Lemma 3.2.2**

*For all $\theta$ and for all $r$, let $M_{\theta,r}$ denote the set of positive-definite, symmetric, d-dimensional matrices $A$ such that $\Theta_r(A) = \theta$. If $(d-1)\pi/2 \leqslant \theta < d\pi/2$ and if $0 < r\tanh(\theta/d) \leqslant 1$, then the function*

$$F_r(A) := \mathrm{Tr}((\mathrm{Id} + r^{-2}A^2)^{-1}(\mathrm{Id} - A^2)) \tag{3.17}$$

*has a unique minimum over $M_{\theta,r}$ at*

$$A_{\theta,r} := r\tanh(\theta/d)\mathrm{Id}. \tag{3.18}$$

*Furthermore, $F_r(A_{\theta,r}) > 0$ unless $r\tan(\theta/d) = 1$.*

**Proof:** Trivially, if the eigenvalues of $A$ are $r\tan(\phi_1), \cdots, r\tan(\phi_d)$, then

$$F_r(A) = \tilde{F}_{r,d}(\phi_1, \cdots, \phi_n) := -dr^2 + \sum_{i=1}^{d}(1+r^2)\cos^2(\phi_i).$$

We therefore denote, for all $m$,

$$E_{\theta,m} := \{(\phi_1, \cdots, \phi_m) \mid \phi_i > 0, \ \phi_1 + \cdots + \phi_m = \theta\},$$

and we use induction on $m$ to show that if $(m-1)\pi/2 \leqslant \theta \leqslant m\pi/2$ and if $0 < r\tan(\theta/m) \leqslant 1$, then $\tilde{F}_{r,m}$ has a unique minimum over $E_{\theta,m}$ at $(\theta/m, \cdots, \theta/m)$.

Suppose first that $m \geqslant 2$ and that $\tilde{F}_{r,m}$ attains its minimum at some point $(\phi_{0,1}, \cdots, \phi_{0,m})$, say, of $E_{m,\theta}$. Using Lagrange multipliers, we verify that there exists $\psi \in [\pi/4, \pi/2]$ such that, for all $i$, $\phi_{0,i} \in \{\psi, \pi/2 - \psi\}$. Let $k$ denote the cardinality of the set of indices $i$ for which $\phi_{0,i} \geqslant \pi/4$. Since $\phi_{0,1} + \cdots + \phi_{0,m} \geqslant (m+1)\pi/2$,

$$k \geqslant (m-1)$$
$$\Rightarrow \quad 2k - m \geqslant 0,$$

and since $\phi_{0,1} + \cdots + \phi_{0,m} = \theta$,

$$k\psi + (m-k)(\pi/2 - \psi) = \theta$$
$$\Leftrightarrow \quad (2k-m)\psi + 2(m-k)\pi/4 = m(\theta/m).$$

Thus, since $\cos^2(\theta)$ is strictly convex over the interval $[\pi/4, \pi/2]$,

$$\tilde{F}_{r,m}(\phi_{0,1}, \cdots, \phi_{0,m}) = (1+r^2)\big((2k-m)\cos^2(\psi) + 2(m-k)\cos^2(\pi/4)\big) - mr^2$$
$$\geqslant (1+r^2)m\cos^2(\theta/m) - mr^2$$
$$= \tilde{F}_{r,m}(\theta/m, \cdots, \theta/m).$$

with equality if and only if $\phi_{0,i} = \theta/m$ for all $i$.

It remains now to verify that $\tilde{F}_{r,m}$ does not attain a strictly lower minimum at any point of $\partial E_{\theta,m}$. For this, it suffices to consider the case where $\phi_m \in \{0, \pi/2\}$. The case $\phi_m = 0$ may only occur when $\theta = (m-1)\pi/2$ and $\phi_i = \pi/2$ for $1 \leqslant i \leqslant (m-1)$. However,

$$\tilde{F}_{r,m}(\pi/2, \cdots, \pi/2, 0) - \tilde{F}_{r,m}(\theta/m, \cdots, \theta/m) = (1+r^2)\cos^2(\theta/m)\big((m-1)\tan^2(\theta/m) - 1\big) \geqslant 0,$$

as desired. When $\phi_m = \pi/2$, by the inductive hypothesis, we have

$$\tilde{F}_{r,m}(\phi_1, \cdots, \phi_{m-1}, \pi/2) - \tilde{F}_{r,m}(\theta/m, \cdots, \theta/m)$$
$$= \tilde{F}_{r,m-1}(\phi_1, \cdots, \phi_{m-1}) - r^2 - \tilde{F}_{r,m}(\theta/m, \cdots, \theta/m)$$
$$\geqslant \tilde{F}_{r,m-1}((\theta - \pi/2)/(m-1), \cdots, (\theta - \pi/2)/(m-1)) - r^2 - \tilde{F}_{r,m}(\theta/m, \cdots, \theta/m)$$
$$= (1+r^2)\big((m-1)\sin^2(\theta'/(m-1)) - m\sin^2(\theta'/m)\big),$$

where $\theta' := m\pi/2 - \theta \in ]0, \pi/2]$, and we leave the reader to verify that this quantity is also always non-negative. We have thus shown that $\tilde{F}_{r,m}$ indeed attains a unique minimum over $E_{\theta,m}$ at $(\theta/m, \cdots, \theta/m)$. It remains only to show that this minimum is non-negative. However, by hypothesis,

$$\tilde{F}_{r,m}(\theta/m, \cdots, \theta/m) = m\cos^2(\theta/m)\big(1 - r^2\tan^2(\theta/m)\big) \geqslant 0,$$

and this completes the proof. $\square$

**Lemma 3.2.3**

If $(d-1)\pi/2 \leqslant \theta_0 < d\pi/2$ and $0 < r\tanh(\theta/d) \leqslant 1$ and if $A : [0, \infty[ \to \Lambda^d$ satisfies

$$\Theta_r(A(0)) = \theta_0, \text{ and} \tag{3.19}$$
$$\dot{A} + A^2 \geqslant \text{Id},$$

then, for all $t \geqslant 0$,

$$\Theta_r(A(t)) \geqslant \theta(t; r, t_0), \tag{3.20}$$

where $\theta(t; r, t_0)$ is as in (3.15).

**Proof:** First, denote

$$\phi(t) := \arctan(\tanh(t_0 + t)/r).$$

Observe that, for all $t$,

$$r\tan(\phi(t)) = \tanh(t + t_0) < 1.$$

Furthermore, this function solves the ordinary differential equation

$$\dot{\phi} = \frac{1}{r}\cos^2(\phi) - r\sin^2(\phi).$$

Now denote $\theta(t) := (1/d)\Theta_r(A(t))$ and denote

$$B(t) := \dot{A}(t) + A(t)^2 \geqslant \text{Id}.$$

Bearing in mind Lemmas 3.1.1 and 3.2.2, for $r\tan(\theta(t)) \leqslant 1$,

$$\begin{aligned}\dot{\theta}(t) &= \frac{1}{rd}\text{Tr}\big((\text{Id} + r^{-2}A(t)^2)^{-1}(B(t) - A(t)^2)\big) \\ &\geqslant \frac{1}{rd}F_r(A(t)) \\ &\geqslant \frac{1}{rd}F_r(r\tan(\theta(t))\text{Id}) \\ &= \frac{1}{r}\cos^2(\theta(t)) - r\sin^2(\theta(t)).\end{aligned}$$

It follows by the classical theory of ordinary differential equations that, for all $t$,

$$\Theta_r(A(t)) = d\theta(t) \geq d\phi(t),$$

as desired. $\square$

**Proof of Theorem 3.2.1:** Indeed, let $\nu : \partial X \to TX$ denote the inward pointing unit normal vector field over $\partial X$. For all $t > 0$, define $e_t : \partial X \to X$ by

$$e_t(x) := \text{Exp}(tx),$$

and observe that $e_t$ is a parametrisation of the level $X_t$ of $X$. For all $t$, let $\nu_t$ denote the upward pointing unit normal vector field over $e_t$ and let $A_t$ denote its corresponding Weingarten operator. Denoting by $R$ the Riemann curvature tensor of $X$, for all $t$, we have (see, for example, Section 2 of [13], under the subheading, "tube formula"),

$$\dot{A}_t = W - A_t^2,$$

where, for all $\xi$,

$$De_t \cdot W\xi = R(\nu_t, De_t \cdot \xi)\nu_t.$$

Since the sectional curvature of $X$ is bounded above by $-1$, $W \geq \text{Id}$, and the result now follows by Lemma 3.2.3. $\square$

**3.3 - Uniqueness.** The estimates of the preceding section yield a preliminary uniqueness result for solutions of the asymptotic Plateau problem. We first require the following estimate.

**Lemma 3.3.1**

Let $X$ be a Cartan–Hadamard manifold with sectional curvature pinched between $-K$ and $-1$. Let $(Y, \phi)$ be a developed Cartan–Hadamard end in $X$ with height function $h$. For $k < 1$, if $Z$ is a graph of constant special Lagrangian curvature equal to $k$ which meets every vertical line then, with $t_0$ as in Lemma 2.1.1,

$$\sup_{x \in Z} h(x) \leqslant \operatorname{arctanh}(k) + t_0. \tag{3.21}$$

**Proof:** Let $\gamma : ]0, \infty[ \to Y$ be a height-parametrised vertical line in $Y$ and let $\Gamma$ denote its image. For all $t > 0$, let $(\mathrm{HS}_t, \alpha_t)$ denote the open half-space in $Y$ centred on $\gamma(t)$. Choose $\epsilon > 0$, and denote $\mathrm{CE} := \mathrm{CE}(\gamma(\epsilon))$. Recall that CE is a Cartan–Hadamard end with height function given by distance to the boundary. Suppose first that $Z \cap \alpha(\mathrm{CE})$ is empty. Then, by Lemma 2.1.1,

$$Z \cap \alpha(\mathrm{HS}_{\gamma(\epsilon)+t_0}) \subseteq Z \cap \alpha(\mathrm{CE}) = \emptyset,$$

so that

$$Z \cap \Gamma \subseteq \gamma(]0, \epsilon + t_0]) \subseteq h^{-1}(]0, \epsilon + t_0]).$$

Now suppose that $Z \cap \alpha(\mathrm{CE})$ is non-empty. For all $s$, let $\mathrm{CE}_s$ denote the level of CE at height $s$. Since $Z$ meets every vertical line, by Lemma 2.7.1, $Z \cap \alpha(\mathrm{CE})$ is relatively compact in $Z$. There therefore exists $s_0$ supremal such that $\alpha(\mathrm{CE}_{s_0})$ meets $Z$. In particular, $\alpha(\mathrm{CE}_{s_0})$ is an exterior tangent to $Z$ at some point. Furthermore, by Lemma 2.1.1 again,

$$Z \cap \alpha(\mathrm{HS}_{\epsilon+s_0+t_0}) \subseteq Z \cap \alpha(\mathrm{CE}_{s_0}) = \emptyset,$$

so that

$$Z \cap \Gamma \subseteq \gamma(]0, \epsilon + s_0 + t_0]) \subseteq h^{-1}(]0, \epsilon + s_0 + t_0]).$$

However, by Corollary 2.5.2, $\alpha(\mathrm{CE}_{s_0})$ has special lagrangian curvature bounded below in the weak sense by $\tanh(s_0)$, so that, by the geometric maximum principle,

$$s_0 \leqslant \operatorname{arctanh}(k).$$

It follows that

$$Z \cap \Gamma \subseteq h^{-1}(]0, \epsilon + \operatorname{arctanh}(k) + t_0]),$$

and the result now follows upon letting $\epsilon$ tend to zero. $\square$

**Lemma 3.3.2**

Let $X$ be a Cartan–Hadamard manifold with sectional curvature pinched between $-K$ and $-1$. Choose $0 < k < 1$ and let $(Y, \phi)$ be a developed Cartan–Hadamard end in $X$ with smooth boundary of constant special lagrangian curvature equal to $k$. If $Z \subseteq Y$ is a graph also of constant special lagrangian curvature equal to $k$ which meets every vertical line, then $Z = \partial Y$.

**Proof:** Indeed, let $h$ denote the height function of $Y$. By Lemma 3.3.1,

$$C := \sup_{x \in \Sigma} h(x) < \infty.$$

Suppose now that $C > 0$. Let $\tilde{\theta} : Y \to \mathbb{R}$ be such that, for all $y$, $\theta(y)$ is the $r$-special lagrangian angle of the level of $Y$ at $y$, where

$$r := k/\tan((d-1)\pi/2d).$$

On the asymptotic Plateau problem in Cartan-Hadamard manifolds.

By Theorem 3.2.1, there exists $\epsilon > 0$ such that, over $h^{-1}([C/2, \infty[)$

$$\tilde{\theta} \geqslant \frac{(d-1)\pi}{2} + \epsilon.$$

Let $\nu$ and $A$ denote respectively the upward pointing unit normal and the corresponding Weingarten operator of $Z$. Let $\nabla$ and Hess denote respectively the gradient and hessian operators of $Z$. Likewise, let $\overline{\nabla}$ and $\overline{\text{Hess}}$ denote respectively the gradient and hessian operators of $Y$. By Lemma 2.5.1 and the subsequent remark, there exists $B \geqslant 0$ such that, over $h^{-1}([C/2, \infty[)$,

$$\|\overline{\text{Hess}}(h)\| \leqslant B. \tag{3.22}$$

Observe now that $Z \cap h^{-1}([C/2, \infty[)$ is complete and, by convexity, has sectional curvature bounded below. It follows by Omori's maximum principal (see Appendix B) that, for all $\delta > 0$, there exists $z_\delta \in Z$ such that

$$\begin{aligned} h(z_\delta) &> C - \delta, \\ \|\nabla h(z_\delta)\| &< \delta, \text{ and} \\ \text{Hess}(h)(z_\delta) &< \delta \text{Id}. \end{aligned} \tag{3.23}$$

In particular, since

$$\langle \nu, \overline{\nabla} h \rangle^2 = 1 - \|\nabla h\|^2, \tag{3.24}$$

for sufficiently small $\delta$,

$$\langle \nu(z_\delta), \overline{\nabla} h(z_\delta) \rangle^{-1} \text{Tr}((1/r)(\text{Id} + r^{-2} A(z_\delta)^2)^{-1} \text{Hess}(h)(z_\delta)) < \frac{\epsilon}{2}. \tag{3.25}$$

By (3.22), (3.23), (3.24) and (B.1) for sufficiently small $\delta$, we have,

$$\Theta_r(\langle \nu(z_\delta), \overline{\nabla} h(z_\delta) \rangle^{-1} \overline{\text{Hess}}(h)(z_\delta)|_{TZ}) \geqslant \tilde{\theta}(z_\delta) - \frac{\epsilon}{2} \geqslant \frac{(d-1)\pi}{2} + \frac{\epsilon}{2}.$$

Thus, by concavity,

$$\begin{aligned} &\langle \nu(z_\delta), \overline{\nabla} h(z_\delta) \rangle^{-1} \text{Tr}((1/r)(\text{Id} + r^{-2} A(z_\delta)^2)^{-1} \text{Hess}(h)(z_\delta)) \\ &= \text{Tr}\big((1/r)(\text{Id} + r^{-2} A(z_\delta)^2)^{-1} \big(\langle \nu(z_\delta), \overline{\nabla} h(z_\delta) \rangle^{-1} \overline{\text{Hess}}(h)(z_\delta)|_{TZ} - A(z_\delta)\big)\big) \\ &\geqslant \Theta_r(\langle \nu(z_\delta), \overline{\nabla} h(z_\delta) \rangle^{-1} \overline{\text{Hess}}(h)(z_\delta)|_{TZ}) - \Theta_r(A(z_\delta)) \\ &\geqslant \frac{\epsilon}{2}. \end{aligned}$$

This contradicts (3.25), and it follows that $C = 0$, as desired. $\square$

It will be helpful to reformulate Lemma 3.3.2 in terms of immersions. Thus, given two quasicomplete ISC immersions $(N_1, e_1)$ and $(N_2, e_2)$ in $M$, we say that $(N_1, e_1)$ **dominates** $(N_2, e_2)$ whenever there exists an open subset $\Omega \subseteq N_1$, a function $\alpha : \Omega \to N_2$ which is a diffeomorphism onto its image, and a smooth function $f : \Omega \to [0, \infty[$ such that

$$e_2 \circ \alpha = \text{Exp}(f\nu_1), \tag{3.26}$$

where $\nu_1$ here denotes the outward-pointing unit normal vector field over $e_1$. Trivially, if we denote by $(\hat{N}_1, \hat{e}_1)$ the end of $(N_1, e_1)$, then this holds if and only if there exists a smooth embedding $\tilde{e}_2 : N_2 \to \hat{N}_1$, whose image is a graph over $\Omega$, such that

$$e_2 = \hat{e}_1 \circ \tilde{e}_2. \tag{3.27}$$

In other words, $(N_1, e_1)$ dominates $(N_2, e_2)$ if and only if $(N_2, e_2)$ factors through the end of $(N_1, e_1)$. Lemma 3.3.2 is now reformulated as follows

**Lemma 3.3.3**

Let $X$ be a Cartan–Hadamard manifold with sectional curvature pinched between $-K$ and $-1$ and let $(Y, \phi)$ be an asymptotic Plateau problem in $\partial_\infty X$. For $k \in ]0, 1[$, if $e_1, e_2 : Y \to X$ are quasicomplete ISC immersions of constant special lagrangian curvature equal to $k$ which solve $(Y, \phi)$, and if $e_1$ dominates $e_2$, then $e_1$ equals $e_2$.

**3.4 - Towards an inverse of the ideal boundary operator.** In the light of Lemma 3.3.2, it will be useful to have a means of determining when solutions of the asymptotic Plateau problem are dominated by over given solutions. We achieve this by constructing a partial inverse to the ideal boundary operator for Cartan–Hadamard ends constructed in Section 2.6. We first require the following technical result.

**Lemma 3.4.1**

Let $X$ be a Cartan–Hadamard manifold of sectional curvature bounded above by $-1$ and let $(Y, \phi)$ be an asymptotic Plateau problem in $\partial_\infty X$. If $\Omega$ is a subset of $Y$ and if $e : \Omega \to Y$ is a quasicomplete ISC immersion such that $\phi_e = \phi$ then, for all $y \in \partial\Omega$,
$$\lim_{z \to y} e(z) = \phi(y). \tag{3.28}$$

**Proof:** Indeed, let $(\hat{\Omega}, \hat{e})$ denote the end of $(\Omega, e)$, let $(\partial_\infty \hat{\Omega}, \partial_\infty \hat{e})$ denote its ideal boundary, and let $\pi_\infty : \Omega \to \partial_\infty \hat{\Omega}$ denote the projection along vertical lines. Recall that $\pi_\infty$ maps $\Omega$ homeomorphically onto $\partial_\infty \hat{\Omega}$ and
$$\partial_\infty \hat{e} \circ \pi_\infty = \phi_e = \phi.$$

Let $(y_m)_{m \in \mathbb{N}}$ be a sequence in $\Omega$ converging to $y$ such that $(e(y_m))_{m \in \mathbb{N}}$ does not converge to $x := \phi(y)$. Choose $\epsilon > 0$, for all $m$, let $(\text{AB}_m, \alpha_m)$ denote the asymptotic ball of $(\partial_\infty \hat{\Omega}, \partial_\infty \hat{e})$ centred on $(y_m, \epsilon)$, and observe that $(\text{AB}_m, \pi_\infty^{-1} \circ \alpha_m)$ is an asymptotic ball of $(\Omega, \phi)$ about $y_m$. By compactness, we may suppose that the sequence $(\hat{e}(y_m, \epsilon))_{m \in \mathbb{N}}$ converges to some point $x'$, say, of $X \cup \partial_\infty X$ distinct from $x$. Since $(\phi(y_m))_{m \in \mathbb{N}}$ converges to $\phi(y) = x$, there thus exists a neighbourhood $V$ of $x$ such that, for sufficiently large $m$, $V \subseteq \text{AB}_m$. Upon reducing $V$ if necessary, we may suppose that it is the homeomorphic image through $\phi$ of some connected neighbourhood $U$ of $y$. Observe, however, that, for all $m$, $(\pi_\infty^{-1} \circ \alpha_m)(\text{AB}_m)$ is the connected component $\phi^{-1}(\text{AB}_m)$ containing $y_m$. Consequently, when $y_m \in U$,
$$y \in U \subseteq (\pi_\infty^{-1} \circ \alpha_m)(\text{AB}_m) \subseteq \pi_\infty^{-1}(\partial_\infty \hat{\Omega}) = \Omega,$$

which is absurd, and the result follows. $\square$

**Lemma 3.4.2**

Let $X$ be a Cartan–Hadamard manifold with sectional curvature pinched between $-K$ and $-1$. Let $(Y, \phi)$ and $(Y', \phi')$ be developed Cartan–Hadamard ends in $X$. Choose $k > 0$ and suppose that $\partial Y$ has smooth boundary of constant $\theta$-special Lagrangian curvature equal to $k$ and that every level of $Y'$ has $\theta$-special Lagrangian curvature strictly greater than $k$ in the weak sense. If $\alpha : \partial_\infty Y' \to \partial_\infty Y$ is an injective local homeomorphism such that
$$\partial_\infty \phi' = \partial_\infty \phi \circ \alpha, \tag{3.29}$$

and if $\alpha(\partial_\infty Y')$ is relatively compact as a subset of $\partial_\infty Y$, then there exists a unique injective local isometry $\tilde{\alpha} : Y' \to Y$ such that
$$\begin{aligned} \phi' &= \phi \circ \tilde{\alpha}, \text{ and} \\ \partial_\infty \tilde{\alpha} &= \alpha. \end{aligned} \tag{3.30}$$

**Proof:** Indeed, let $\gamma : ]0, \infty[ \to Y'$ be a height-parametrised vertical line with end point $\gamma_\infty$. For sufficiently large $t > 0$, there exists a unique lift $\tilde{\gamma} : [t, \infty[ \to Y$ with end point $\tilde{\gamma}_\infty$ such that

$$\begin{aligned} \phi \circ \tilde{\gamma} &= \phi' \circ \gamma, \text{ and} \\ \tilde{\gamma}_\infty &= \alpha(\gamma_\infty). \end{aligned}$$

Since $\alpha(\partial_\infty Y')$ is relatively compact in $\partial_\infty Y$, it follows by Lemma 3.4.1 that this lift can be performed uniformly over the leaf space of $Y'$, even near the boundary. Thus, denoting by $h$ the height function of $Y'$, there exists $t > 0$ and a unique function $\tilde{\alpha} : h^{-1}([t, \infty[) \to Y$ such that

$$\phi \circ \tilde{\alpha} = \phi', \text{ and}$$
$$\partial_\infty \tilde{\alpha} = \alpha.$$

Now let $t_0$ be infimal such that $\alpha$ is defined over $h^{-1}([t_0, \infty[)$ and suppose that $t_0 > 0$. For all $t$, let $Y'_t$ denote the level of $Y'$ at height $t$. By Lemma 2.4.6, for all $t > t_0$, $\alpha(Y'_t)$ is a graph over $\Omega := \alpha(\partial_\infty Y'_t)$. From this we see that the only obstruction to $\alpha$ being extended beyond $t_0$ is for $\alpha(X'_{t_0})$ to be an exterior tangent to $\partial X$ at some point. Since this is impossible by the geometric maximum principle, the result follows. □

## 4 - Compactness and existence.

**4.1 - Cheeger–Gromov convergence.** In [35], building on the work [24] of Labourie, we proved a general compactness result for families of quasicomplete hypersurfaces of constant special lagrangian curvature immersed in riemannian manifolds. We expressed these results using the concept of Cheeger–Gromov convergence which we now review (c.f. [29] and [33]). First, we define a **pointed riemannian manifold** to be a triplet $(X, g, x)$, where $X$ is a smooth manifold, $g$ is a riemannian metric, and $x$ is a point of $X$. We say that a sequence $(X_m, g_m, x_m)_{m \in \mathbb{N}}$ of complete pointed riemannian manifolds converges to the complete pointed riemannian manifold $(X_\infty, g_\infty, x_\infty)$ in the **Cheeger–Gromov sense** whenever there exists a sequence $(\Phi_m)_{m \in \mathbb{N}}$ of functions such that

(1) for all $m$, $\Phi_m : X_\infty \to X_m$ and $\Phi_\infty(x_\infty) = x_m$; and

for every relatively compact open subset $\Omega$ of $X_\infty$, there exists $M$ such that

(2) for all $m \geq M$, the restriction of $\Phi_m$ to $\Omega$ defines a smooth diffeomorphism onto its image; and

(3) the sequence $((\Phi_m|_\Omega)^* g_m)_{m \geq M}$ converges to $g_\infty|_\Omega$ in the $C^\infty_{\text{loc}}$ sense.

We call $(\Phi_m)_{m \in \mathbb{N}}$ a sequence of **convergence maps** of $(X_m, g_m, x_m)_{m \in \mathbb{N}}$ with respect to $(X_\infty, g_\infty, x_\infty)$. Although the convergence maps are trivially non-unique, it is straightforward to show that any two such sequences $(\Phi_m)_{m \in \mathbb{N}}$ and $(\Phi'_m)_{m \in \mathbb{N}}$ are equivalent in the sense that there exists an isometry $\Psi : X_\infty \to X_\infty$ preserving $x_\infty$ such that, for any two relatively compact open subsets $U \subseteq \overline{U} \subseteq V$ of $X_\infty$, there exists $M$ such that

(1) for all $m \geq M$, the respective restrictions of $\Phi_m$ and $\Phi'_m \circ \Psi$ to $U$ and $V$ define smooth diffeomorphisms onto their images;

(2) for all $m \geq M$, $(\Phi'_m \circ \Psi)(U) \subseteq \Phi_m(V)$; and

(3) the sequence $((\Phi_m|_V)^{-1} \circ \Phi'_m \circ \Psi)_{m \geq M}$ converges in the $C^\infty$ sense to the identity map over $U$.

The concept of Cheeger–Gromov convergence applies to sequences of immersed submanifolds as follows. We say that a sequence $(Y_m, x_m, \phi_m)_{m \in \mathbb{N}}$ of complete pointed immersed submanifolds in a complete riemannian manifold $(X, g)$ converges to the complete pointed immersed submanifold $(Y_\infty, x_\infty, \phi_\infty)$ in the **Cheeger–Gromov sense** whenever $(Y_m, x_m, \phi_m^* g)_{m \in \mathbb{N}}$ converges to $(Y_\infty, x_\infty, \phi_\infty^* g)$ in the Cheeger–Gromov sense and for one, and thus for any, sequence $(\Phi_m)_{m \in \mathbb{N}}$ of convergence maps, the sequence $(\phi_m \circ \Phi_m)_{m \in \mathbb{N}}$ converges to $\phi_\infty$ in the $C^\infty_{\text{loc}}$ sense.

These definitions are readily extended in a number of ways. For example, in the case of a sequence $(X_m, g_m, x_m)_{m \in \mathbb{N}}$ of pointed riemannian manifolds, the hypothesis of completeness is unnecessary. Instead, it is sufficient to assume that for all $R > 0$, there exists $M$ such that, for all $m \geq M$, the closed ball of radius $R$ about $x_m$ in $(X_m, g_m)$ is compact. Likewise, in the case of immersed submanifolds, the target space can be replaced with a sequence $(X_m, g_m, x_m)_{m \in \mathbb{N}}$ of pointed riemannian manifolds converging in the Cheeger–Gromov sense to some complete limit. Furthermore, it is not necessary to suppose that the riemannian manifolds in this sequence are complete, and so on.

Following [25], given an ISC immersed hypersurface $(Y, e)$ in a riemannian manifold $X$, we define its **Gauss lift** $\hat{e} : Y \to SX$ by

$$\hat{e} := \nu_e, \tag{4.1}$$

where $\nu_e$ here denotes the unit normal vector field of $e$ compatible with the orientation. We use this notation in order to emphasize our interest in this function, not as a vector field over $e$, but as an immersion in its own right.

**Theorem 4.1.1**

Let $X$ be a complete riemannian manifold. Choose $k > 0$ and let $(Y_m, e_m, y_m)_{m \in \mathbb{N}}$ be a sequence of quasicomplete, pointed, ISC immersed hypersurfaces in $X$ of constant special lagrangian curvature equal to $k$ and, for all $m$, let $\hat{e}_m : Y_m \to SX$ denote the Gauss lift of $e_m$. If $(\hat{e}_m(y_m))_{m \in \mathbb{N}}$ remains within a compact subset of $SX$, then there exists a complete pointed immersed submanifold $(Y_\infty, \hat{e}_\infty, y_\infty)$ in $SX$ towards which $(Y_m, \hat{e}_m, y_m)_{m \in \mathbb{N}}$ subconverges in the Cheeger–Gromov sense.

**Remark 4.1.1.** We have stated Theorem 4.1.1 here in its simplest possible form. However, all the generalisations of Cheeger–Gromov convergence discussed above also apply here. We may also suppose that, for all $m$, the special lagrangian curvature of $(Y_m, e_m)$ is equal to $k_m$, where $(k_m)_{m \in \mathbb{N}}$ converges to some non-zero limit $k_\infty$, and so on.

Significantly, Theorem 4.1.1 does not state that the limit submanifold $(\hat{Y}_\infty, \hat{e}_\infty)$ is itself a lift of some quasicomplete immersed hypersurface. However, a remarkable result [24] of Labourie states that there exists only one, very precise, mode of degeneration. Our extension to the higher dimensional case, proven in [35], forms the content of the next theorem. First, given a complete geodesic $\Gamma$ in $X$, we denote by $S\Gamma \subseteq SX$ the set of unit normal vectors over $\Gamma$. Observe that $S\Gamma$ is an immersed submanifold isometric to $\mathbb{S}^{d-1} \times \mathbb{R}$. We now define a **tube** to be any immersed submanifold $(Y, \hat{e})$ in $SX$ which is a cover of $S\Gamma$, for some complete geodesic $\Gamma$.

**Theorem 4.1.2**

If $(Y_\infty, \hat{e}_\infty)$ is a limit of a sequence of lifts of quasicomplete ISC immersed hypersurfaces of constant special lagrangian curvature as in Theorem 4.1.1, then either $(Y_\infty, \hat{e}_\infty)$ is a tube or $(Y_\infty, \pi \circ \hat{e}_\infty)$ is a quasicomplete ISC immersed hypersurface of constant special lagrangian curvature.

**Remark 4.1.2.** In [24] Labourie calls such hypersurfaces "surfaces rideaux", that is "curtain surfaces".

The phenomenon described in Theorem 4.1.2 is illustrated by the case of equidistant cylinders about some geodesic $\Gamma$ in hyperbolic space. Indeed, for all $r$, let $C_r$ denote the cylinder of points in $\mathbb{H}^{d+1}$ lying at distance $r$ from $\Gamma$. For all $r$, $C_r$ has constant special Lagrangian curvature equal to $k(r)$, for some function $k$ of $r$. In fact, it is a nice exercise to show that

$$k(r) = \frac{1}{\sqrt{d-1}} \tan((d-1)\pi/2d) + \mathrm{O}(r^2). \tag{4.2}$$

The family $(C_r)_{r>0}$ trivially degenerates as $r$ tends to zero. On the other hand, the family $(\hat{C}_r)_{r>0}$ of lifts of these cylinders converges towards $S\Gamma$, in accordance with Theorems 4.1.1 and 4.1.2. Furthermore, in the 2-dimensional case, for all $m$, the family $(\hat{C}_{m,r})_{r>0}$ of $m$-fold covers of these lifts converges towards an $m$-fold cover of this tube. Finally, the sequence $(\hat{C}_{m,1/m})_{m \in \mathbb{N}}$ converges towards the universal cover of $S\Gamma$. Theorem 4.1.2 affirms that these are the only modes of degeneration that can occur.

**4.2 - Monotone convergence.** We are now ready to prove one of the main components of this paper, namely a monotone convergence result for sequences of solutions to the asymptotic Plateau problem.

**Theorem 4.2.1, Monotone convergence**

Let $X$ be a Cartan–Hadamard manifold of sectional curvature pinched between $-K$ and $-1$. Let $(Y, \phi)$ be an asymptotic Plateau problem in $\partial_\infty X$ which is not of finite type 0, 1 or 2. Let $(\Omega_m)_{m \in \mathbb{N}}$ be a nested sequence of relatively compact open subsets of $Y$ which exhausts $Y$. If, for $k \in ]0, 1[$ and for all $m$, there exists a quasicomplete ISC immersion $e_m : \Omega_m \to X$ of constant special lagrangian curvature equal to $k$ solving $(\Omega_m, \phi)$, then $(e_m)_{m \in \mathbb{N}}$ subconverges in the $C^0_{\text{loc}}$ sense to a quasicomplete ISC immersion $e : Y \to X$ of constant special lagrangian curvature equal to $k$ solving $(Y, \phi)$.

**Remark 4.2.1.** For ease of presentation, Theorem 4.2.1 is stated and proven in its simplest form. More generally, it is not necessary to suppose quasicompleteness of $e_m$ for all $m$. Indeed, it is sufficient that, for all $y \in Y$, for all $R > 0$, and for all sufficiently large $m$, the closed ball of radius $R$ about $y$ in $\Omega_m$ with respect to the metric $\text{I}_m + \text{III}_m$ is compact, where here $\text{I}_m$ and $\text{III}_m$ denote respectively the first and third fundamental forms of $e_m$.

The proof of Theorem 4.2.1 is given in the following series of lemmas.

**Lemma 4.2.2**

For all $y \in Y$, and for every asymptotic ball $(AB, \alpha)$ in $Y$ containing $y$, there exists a neighbourhood $U$ of $\phi(y)$ in $X \cup \partial_\infty X$ such that, for all sufficiently large $m$,

$$e_m(y) \notin U. \tag{4.3}$$

**Proof:** Let $\xi \in SX$ be such that $AB = AB(\xi)$ and, for all $r > 0$, denote

$$\text{CE}_r := \{x \in \text{CE}(\xi) \mid d(x, \partial \text{CE}(\xi)) > r\}.$$

Let $U$ be a neighbourhood of $\phi(y)$ in $X \cup \partial_\infty X$ with closure contained in $\text{CE}_{r_0}$, where

$$r_0 := \operatorname{arctanh}(k).$$

We show that $U$ has the desired properties. Indeed, without loss of generality, we may suppose that, for all $m$, $\alpha(AB) \subseteq \Omega_m$. For all $m$, let $(\hat{\Omega}_m, \hat{e}_m)$ denote the end of $(\Omega_m, e_m)$. Let $\pi_\infty : \hat{\Omega}_m \to \partial_\infty \hat{\Omega}_m$ denote the projection along vertical lines and observe that, by definition,

$$\partial_\infty \hat{e}_m \circ (\pi_\infty \circ \alpha) = (\partial_\infty \hat{e}_m \circ \pi_\infty) \circ \alpha = \phi \circ \alpha = \text{Id}.$$

By Corollary 2.5.2, every level of $\text{CE}_{r_0}$ has special lagrangian curvature strictly greater than $k$ in the weak sense. Since $\partial_\infty \text{CE}_{r_0} = AB$, it follows by Lemma 3.4.2 that there exists a local isometry $\hat{\alpha} : \text{CE}_{r_0} \to \hat{\Omega}_m$ such that

$$\hat{e}_m \circ \hat{\alpha} = \text{Id}, \text{ and}$$
$$\partial_\infty \hat{\alpha} = \pi_\infty \circ \alpha.$$

By Lemma 2.4.6, for all $\epsilon > 0$, the level $\hat{\alpha}(\text{CE}_{r_0, \epsilon})$ is a graph in $\hat{\Omega}_m$. In particular, the vertical line of $\hat{\Omega}_m$ over $y$ intersects this level at some point. In other words, the geodesic ray $\gamma : [0, \infty[ \to X; t \mapsto \hat{e}_m(y, t)$ crosses $\text{CE}_{r_0, \epsilon} = \partial \text{CE}_{r_0, \epsilon}$ at some point. However, for sufficiently small $\epsilon$,

$$\partial \text{CE}_{r_0 + \epsilon} \cap U = \emptyset,$$

so that, by convexity,

$$e_m(y) = \hat{e}_m(y, 0) = \gamma(0) \notin U,$$

as desired. $\square$

**Lemma 4.2.3**

For all $y \in Y$, there exists a compact subset $K$ of $X$ such that, for all $m$ with $y \in \Omega_m$,

$$e_m(y) \in K. \tag{4.4}$$

**Proof:** Choose $y \in Y$ and let $U$ be a neighbourhood of $y$, homeomorphic to a ball, over which $\phi$ restricts to a homeomorphism onto its image. Let $(AB, \alpha)$ be an asymptotic ball in $Y$ with relatively compact image in $U$. We may suppose that, for all $m$, $\alpha(AB)$ is also a relatively compact subset of $\Omega_m$. Let $\omega$ denote the Kulkarni–Pinkall section of $\phi$ and, for all $m$, let $\omega_m$ denote the Kulkarni–Pinkall section of the restriction of $\phi$ to $\Omega_m$. By Lemma 2.3.1,

$$\omega_m \leqslant \omega.$$

Thus, by Theorem 2.8.2, for all $m$,

$$e_m(y) \in \overline{\omega}_m(y) \subseteq \overline{\omega}(y).$$

Consequently, for all $m$,

$$e_m(y) \in K := \overline{\omega}(y) \setminus V,$$

where $V$ is as in Lemma 4.2.2. Finally, since $(Y, \phi)$ is not of finite type $0$ or $1$, $K$ is compact, and the result follows. $\square$

Now fix $y_0 \in Y$. Without loss of generality, we may suppose that $y_0 \in \Omega_m$ for all $m$. For all $m$, let $\hat{e}_m : \Omega_m \to SX$ denote the Gauss lift of $e_m$. By Lemma 4.2.3 and Theorem 4.1.1, we may suppose that there exists a complete pointed immersed submanifold $(Y', \hat{e}', y')$ towards which $(\Omega_m, \hat{e}_m, y_0)_{m \in \mathbb{N}}$ converges in the Cheeger–Gromov sense. Let $(\beta_m)_{m \in \mathbb{N}}$ be a sequence of convergence maps of $(\Omega_m, \hat{e}_m, y_0)$ with respect to this limit.

**Lemma 4.2.4**

For every compact subset $K$ of $Y'$, there exists a compact subset $L$ of $Y$ such that, for all $m$,

$$\beta_m(K) \subseteq L. \tag{4.5}$$

**Proof:** First, let $\phi' : Y' \to \partial_\infty X$ denote the asymptotic Gauss map of $\hat{e}'$. That is, $\phi' := \text{Hor} \circ \hat{e}'$, where $\text{Hor} : X \to \partial_\infty X$ here denotes the horizon map of $X$. Trivially, we may suppose that $K$ is connected and contains $y'$. Let $(AB'_i, \alpha'_i)_{1 \leq i \leq k}$ be a finite set of asymptotic balls in $(Y', \phi')$ covering $K$ and indexed such that $y' \in \alpha'_1(AB_i)$ and that, for all $i > 1$,

$$\alpha'_i(AB_i) \cap \left( \bigcup_{1 \leqslant j < i} \alpha'_j(AB_j) \right) \neq \emptyset.$$

Suppose furthermore that, for all $i$, the closure of $(AB_i, \alpha'_i)$ is contained in the strictly larger ball $(\widetilde{AB}_i, \tilde{\alpha}'_i)$ whose closure is itself also a compact subset of $Y'$. For all $i$, $(\phi \circ \beta_m \circ \tilde{\alpha}'_i)$ converges locally uniformly to $(\phi' \circ \tilde{\alpha}'_i) = \text{Id}$. Thus, for sufficiently large $m$, this function has an inverse $\gamma_{i,m}$ defined over a neighbourhood of $\overline{AB}_i$. For all such $m$, denote

$$\alpha_{m,i} := (\beta_m \circ \tilde{\alpha}'_i \circ \gamma_{i,m}).$$

In particular, for all such $m$,

$$\phi \circ \alpha_{m,i} = (\phi \circ \beta_m \circ \tilde{\alpha}'_i) \circ \gamma_{i,m} = \text{Id},$$

so that $(AB_i, \alpha_{m,i})$ is an asymptotic ball in $(Y, \phi)$. We now claim that, for sufficiently large $m$, $\alpha_{m,i}$ is independent of $m$. Indeed, by hypothesis, $y' \in \alpha_{m,1}(AB_1)$ for all $m$, so that $\alpha_1 := \alpha_{m,1}$ is independent of $m$. Now suppose that $\alpha_j := \alpha_{m,j}$ is independent of $m$ for all $1 \leqslant j \leqslant i - 1$. By construction, there exists $j < i$ such that

$$\alpha'_j(AB_j) \cap \alpha'_i(AB_i) \neq \emptyset.$$

Let $x \in \mathrm{AB}_j$ be such that $\alpha'_j(x)$ lies in this intersection. In particular, upon composing with $\phi'$, we see that $x \in \mathrm{AB}_i$. Let $V$ be a connected neighbourhood of $x$ with closure contained in $\mathrm{AB}_i \cap \mathrm{AB}_j$ and observe that

$$\alpha'_j|_V = \alpha'_i|_V.$$

Since $(\gamma_{j,m})_{m \in \mathbb{N}}$ and $(\gamma_{i,m})_{m \in \mathbb{N}}$ converge locally uniformly to the identity, we may suppose that, for all sufficiently large $m$,
$$x \in \gamma_{j,m}^{-1}(V), \text{ and}$$
$$V \subseteq \gamma_{i,m}(\mathrm{AB}_i).$$

For all such $m$,

$$\alpha_j(x) = \alpha_{j,m}(x) = (\beta_m \circ \tilde{\alpha}'_j \circ \gamma_{j,m})(x) \subseteq (\beta_m \circ \tilde{\alpha}'_i \circ \gamma_{i,m})(\mathrm{AB}_i) = \alpha_{i,m}(\mathrm{AB}_i),$$

from which it follows that, for all sufficiently large $m$, $\alpha_i := \alpha_{m,i}$ is also independent of $m$, and the assertion follows by induction. In particular, for all sufficiently large $m$,

$$\beta_m(K) \subseteq \bigcup_{i=1}^{k} \overline{\alpha_i(\mathrm{AB}_i)} =: L,$$

as desired. □

**Lemma 4.2.5**

*The sequence $(\beta_m)_{m \in \mathbb{N}}$ subconverges locally uniformly to a homeomorphism $\beta_\infty : Y' \to Y$.*

**Proof:** We first show that $\beta_m$ subconverges to an injective local homeomorphism $\beta_\infty : Y' \to Y$. Indeed, by Lemma 4.2.4, upon extracting a subsequence if necessary, we may suppose that $(\beta_m)_{m \in \mathbb{N}}$ converges over a countable dense subset $Y'_0$ of $Y'$ to some function $\beta_{\infty,0} : Y'_0 \to Y$. Since $(\phi \circ \beta_m)_{m \in \mathbb{N}}$ converges locally uniformly to the local homeomorphism $\phi'$, it then follows that $(\beta_m)_{m \in \mathbb{N}}$ itself converges locally uniformly over the whole of $Y'$ to a local homeomorphism $\beta_\infty : Y' \to Y$. Furthermore, since, for every compact subset $L'$ of $Y'$, and for sufficiently large $m$, the restriction of $\beta_m$ to $L'$ is injective, it follows that $\beta_\infty$ is also injective.

It remains only to show that $\beta_\infty$ is surjective. To this end, denote $\Omega := \beta_\infty(Y')$ and observe that $(\hat{e}_m)_{m \in \mathbb{N}}$ converges locally uniformly over this set to $\hat{e}_\infty := (\hat{e}' \circ \beta_\infty^{-1})$. Suppose now that $\Omega$ is not the whole of $Y$. By connectedness, $\partial \Omega$ is non-empty. Let $y'_0$ be a point of this set and let $\Omega$ and $\hat{e}'_\infty : \Omega' \to SX$ denote the domain and function obtained upon applying the preceding construction about $y'_0$ instead of $y_0$. Since $\hat{e}_\infty$ coincides with $\hat{e}'_\infty$ over $\Omega \cap \Omega'$, it follows that $\hat{e}_\infty$ smoothly extends across a neighbourhood of $y$. This contradicts completeness, and the result follows. □

We now complete the proof of Theorem 4.2.1.

**Proof of Theorem 4.2.1:** It remains only to show that $e_\infty := \pi \circ \hat{e}_\infty$ is an immersion. However, suppose the contrary. By Theorem 4.1.2, $(Y, \hat{e}_\infty)$ is a tube, so that $(Y, \phi)$ is of finite type 2, which is absurd, and this completes the proof. □

**4.3 - Convex cobordisms.** In [36], we solved the classical Plateau problem for hypersurfaces of constant special lagrangian curvature in Cartan–Hadamard manifolds. We will now use this result as the starting point for constructing solutions to the asymptotic Plateau problem. We first recall the concept of convex cobordism, which we introduced in [36], and which provides a useful framework for stating and applying our results. We define a **convex cobordism** in a riemannian manifold $X$ to be a pair $(Y, \phi)$, where $Y$ is a compact riemannian manifold with corners of dimension equal to that of $X$ and $\phi : Y \to X$ is a local isometry such that

(1) the boundary $\partial Y$ of $Y$ consists of two smooth components $\partial_+ Y$ and $\partial_- Y$ which meet transversally along a smooth codimension 2 submanifold;

(2) the **upper boundary** $\partial_+ Y$ is ISC with $Y$ lying on its inside;

(3) the **lower boundary** $\partial_-Y$ is ISC with $Y$ lying on its outside; and

(4) $Y$ is foliated by geodesics normal to its lower boundary.

The structure of a convex cobordism is illustrated in Figure 4.3.3. We call the geodesic segments normal to $\partial_-Y$ **vertical lines**. We call the resulting foliation the **vertical line foliation** of $Y$, and we denote its leaf space by $\mathcal{V}$. We define the **height function** $h$ of $Y$ to be the distance to $\partial_-Y$ along vertical lines. By standard properties of convex subsets of Cartan–Hadamard manifolds, this function is $C^{1,1}_{\text{loc}}$ and strictly convex, and its gradient flow lines are none other than the vertical lines.

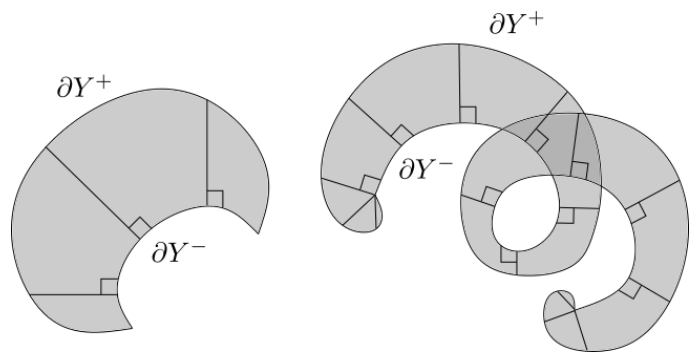

**Figure 4.3.3 - Convex cobordisms -** Observe that the canonical projection along vertical lines $\pi: Y \to \partial Y^-$ does not necessarily define a bijection from $\mathcal{V}$ onto the lower boundary. Indeed, in the second example shown here, multiple vertical lines emanate from boundary points of $\partial_-Y$.

In [36], we prove the following result.

**Theorem 4.3.1**

Let $X$ be a Cartan–Hadamard manifold of sectional curvature bounded above by $-1$ and let $(Y, e)$ be a compact ISC immersed surface with boundary in $X$. Choose $0 < k < 1$ and suppose that $Y$ has special lagrangian curvature at every point strictly greater than $k$, and that the restriction of $e$ to $\partial Y$ does not self-intersect tangentially. Then there exists a unique convex cobordism $(\hat{Y}, \hat{e})$ in $X$ whose upper boundary is a reparametrisation of $(Y, e)$ and whose lower boundary has constant special lagrangian curvature equal to $k$.

We will also require the following degenerate form of the concept of convex cobordisms, which provides a bridge towards the concept of Cartan–Hadamard ends, studied in the preceding chapters. Given a codimension 2 immersed submanifold $Z$ in a riemannian manifold $X$, we say that a unit normal vector $\nu_p$ over $Z$ at some point $p$ is a **convex direction** at this point whenever, for any tangent vector field $\xi$ over $Z$, not vanishing at this point,

$$\langle (\nabla_\xi \xi)(p), \nu_p \rangle < 0. \tag{4.6}$$

The reader may verify that this property holds if and only if $Z$ extends in a neighbourhood of $p$ to an ISC hypersurface with outward-pointing normal $\nu_p$. We now define an **ideal convex cobordism** in a riemannian manifold $X$ to be a pair $(Y, \phi)$ where $Y$ is a riemannian manifold with corners and $\phi : Y \to X$ is a local isometry such that

(1) the boundary $\partial Y$ of $Y$ consists of two smooth components $\partial_0 Y$ and $\partial_-Y$ which meet (transversally) along a smooth codimension 2 submanifold $Z$;

(2) the **lateral boundary** $\partial_0 Y$ is foliated by complete geodesic rays leaving $Z$ in convex normal directions;

(3) the **lower boundary** $\partial_-Y$ is compact and ISC with $Y$ lying on its outside; and

(4) $Y$ is foliated by complete geodesic rays normal to its lower boundary.

This concept is illustrated in Figure 4.3.4. Ideal convex cobordisms have structures similar to those of Cartan–Hadamard ends, which it is again worth reviewing in some detail. First we say that a tangent vector $\xi \in TY$ over some point $y$ is **upward pointing** whenever

$$\langle \xi, \nabla h(y) \rangle > 0.$$

We denote respectively by $T^+Y$ and $S^+Y$ the subbundles of upward pointing tangent vectors and unit upward pointing tangent vectors over $Y$. We call the geodesic rays normal to $\partial_- Y$ **vertical lines**, we call the resulting foliation the **vertical line foliation**, and we denote its leaf space by $\mathcal{V}$. Observe that $\mathcal{V}$ is the union of two disjoint components $\partial \mathcal{V}$ and $\mathcal{V}^o$, where $\partial \mathcal{V}$ foliates the lateral boundary and $\mathcal{V}^o$ foliates the interior of $Y$. We define the **height function** $h$ of $Y$ to be the distance to the lower boundary along vertical lines. By standard properties of convex subsets of Cartan–Hadamard manifolds, this function is $C^{1,1}$ and strictly convex, and its gradient flow lines are none other than the vertical lines. Finally, we call the level sets of $h$ the **levels** of $Y$. These form another foliation of $Y$ by $C^{1,1}$ embedded hypersurfaces all transverse to the lateral boundary which we call the **level set foliation** and which we denote by $(Y_t)_{t \geqslant 0}$. The vertical line foliation and the level set foliation are transverse to one another and every vertical line intersects every level of non-zero height at exactly one point. In particular, every level of non-zero height is naturally homeomorphic to $\mathcal{V}$.

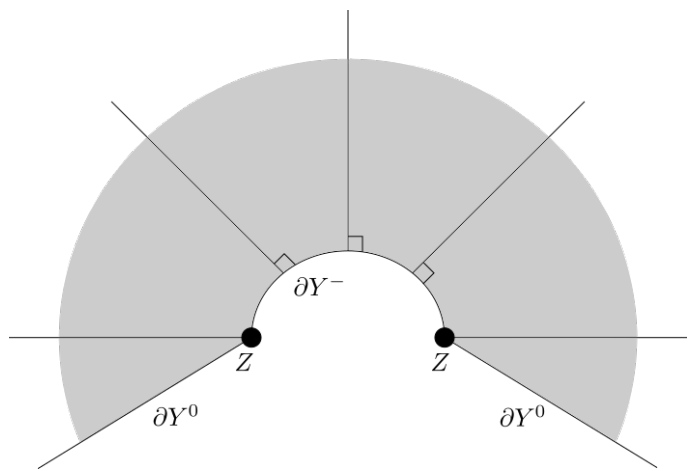

**Figure 4.3.4 - An ideal convex cobordism -** Observe again that multiple vertical lines may emanate from boundary points of $\partial_- Y$.

In the case of ideal convex cobordisms, we have the following restricted version of the Hopf-Rinow theorem.

**Lemma 4.3.2**

Let $X$ be a Cartan–Hadamard manifold with sectional curvature bounded above by $-1$. Let $(Y, \phi)$ be an ideal convex cobordism in $X$ with height function $h$ and let $\mathcal{V}$ denote the leaf space of its vertical line foliation. For every compact subset $K$ of $\mathcal{V}^o$, there exists $h > 0$ such that, for every point $y$ of $Y$ lying on a vertical line of $K$ at height greater than $h_0$, $T_{+,y}Y$ is contained in the domain of the exponential map of $Y$.

**Proof:** Choose $\epsilon > 0$ and let $Y_\epsilon$ denote the level of $Y$ at height $\epsilon$. Let $K_\epsilon \subseteq Y_\epsilon$ denote the set of points lying on vertical lines of $K$. Let $\delta > 0$ be such that

$$d_\epsilon(K_\epsilon, \partial Y_\epsilon) > \delta,$$

where $d_\epsilon$ here denotes the intrinsic distance in $Y_\epsilon$. We verify as in Lemma 2.7.1 that $h_0 > 0$, chosen such that

$$\mathrm{arccosh}(\coth(h_0 - \epsilon)) < \delta,$$

has the desired properties. □

Let $y$ be a point of $Y$ such that $T_{+,y}Y$ is contained in the domain of the exponential map. Denote $x := \phi(y)$ and $\xi := D\phi(y) \cdot \nabla h(y)$ and observe that $\text{Exp}_x \circ D\phi(y)$ maps $T_{+,y}Y$ diffeomorphically onto $\text{HS}(\xi)$. We thus denote

$$\text{HS}(y) := \text{HS}(\xi), \tag{4.7}$$

and we define $\alpha(y) : \text{HS}(y) \to Y$ such that

$$\alpha(y) \circ \text{Exp}_x \circ D\phi(y) = \text{Exp}_y. \tag{4.8}$$

We call $(\text{HS}(y), \alpha(y))$ the **open half-space** of $Y$ centred on $Y$.

We conclude this section by studying the ideal boundary of an ideal convex cobordism. First, as in Lemma 2.6.1, we show that complete, unit speed geodesic rays $\gamma : [0, \infty[ \to Y$ are divided into two classes, namely those that are bounded and those that are unbounded. We thus define the **ideal boundary** $\partial_\infty Y$ of $Y$ to be the space of equivalence classes of unbounded, complete, unit speed geodesic rays in $Y$, where two such rays are deemed equivalent whenever they are asymptotic to one another. We verify as before that every complete, unit speed geodesic ray in $Y$ is asymptotic to a unique vertical line, from which it follows that $\partial_\infty Y$ is naturally homeomorphic to $\mathcal{V}$.

With $y$, $x$ and $\xi$ as before, we now observe that $\text{Hor}_x^X \circ D\phi(y)$ maps $S_{+,y}Y$ homeomorphically onto $\text{AB}(\xi)$. We thus denote

$$\text{AB}(y) := \text{AB}(\xi), \tag{4.9}$$

and we define $\alpha(y) : \text{AB}(y) \to \partial_\infty Y$ such that

$$\alpha(y) \circ \text{Hor}_x^X \circ D\phi(y) = \text{Hor}_y^Y. \tag{4.10}$$

We call $(\text{AB}(y), \alpha(y))$ the **asymptotic ball** in $\partial_\infty Y$ centred on $y$.

**4.4 - Constructing ideal convex cobordisms.** Let $X$ be a Cartan–Hadamard manifold with sectional curvature bounded above by $-1$. A straightforward construction of ideal convex cobordisms in $X$ is given as follows. Let $(Y, e)$ be a compact ISC hypersurface with boundary immersed in $X$ and let $\nu : Y \to SX$ denote its outward pointing unit normal vector field. Define $\hat{Y} := Y \times [0, \infty[$ and $\hat{e} : \hat{Y} \to X$ by

$$\hat{e}(x, t) := \text{Exp}(t\nu(x)).$$

By standard properties of Cartan–Hadamard manifolds, $\hat{e}$ is a local diffeomorphism. We furnish $\hat{Y}$ with the unique riemannian metric making $\hat{e}$ into a local isometry. $(\hat{Y}, \hat{e})$ is trivially an ideal convex cobordism, which we call the ideal convex cobordism **above** $(Y, e)$.

A more sophisticated construction is given as follows. Let $(Y_1, \phi_1)$ and $(Y_2, \phi_2)$ be respectively a convex cobordism and an ideal convex cobordism in $X$. If the upper boundary of $(Y_1, \phi_1)$ is a reparametrisation of the lower boundary of $(Y_2, \phi_2)$, then, as illustrated in Figure 4.4.5, the two join in a natural manner to define a smooth manifold with corners $Y_3$ and a smooth local isometry $\phi_3 : Y_3 \to X$.

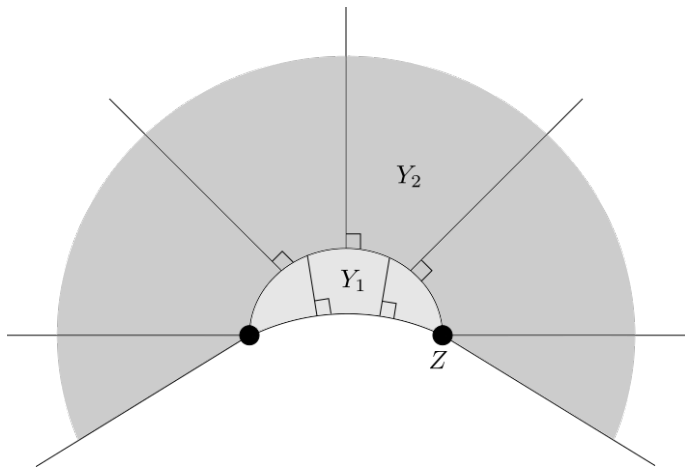

**Figure 4.4.5 - Joining convex cobordisms** - The vertical line foliation of the join coincides with that of $Y_1$ but not with that of $Y_2$.

**Lemma 4.4.1**

Let $X$ be a Cartan–Hadamard manifold. Let $(Y_1, \phi_1)$ and $(Y_2, \phi_2)$ be respectively a convex cobordism and an ideal convex cobordism in $X$. If the upper boundary of $(Y_1, \phi_1)$ is a reparametrisation of the lower boundary of $(Y_2, \phi_2)$, then their natural join $(Y_3, \phi_3)$ is also an ideal convex cobordism in $X$.

**Proof:** The only non-trivial property to be verified is that $Y_3$ is foliated by complete geodesic rays normal to its lower boundary $\partial_- Y_3$. To this end, let $h : Y_3 \to [0, \infty[$ denote the distance in $Y_3$ to $\partial_- Y_3$. We will show that $h$ is a height function and the result will follow with the vertical line foliation given by its gradient flow lines.

By standard properties of convex subsets of Cartan–Hadamard manifolds, it is in fact sufficient to show that every $y \in Y_3$ has a unique closest point in $\partial_- Y_3$. We therefore suppose that contrary. Observe first that, by convexity and compactness of $\partial_- Y_3$, there exists $\delta > 0$ such that, for any two geodesic segments $\gamma_1, \gamma_2 : [0, R] \to Y_3$ normal to $\partial_- Y_3$ at $t = 0$ and such that $d(\gamma_1(0), \gamma_2(0)) < \delta$,

$$(\phi_3 \circ \gamma_1)([0, R]) \cap (\phi_3 \circ \gamma_2)([0, R]) = \emptyset$$
$$\Rightarrow \gamma_1([0, R]) \cap \gamma_2([0, R]) = \emptyset.$$

Consequently, if any point $y \in Y_3$ has two distinct closest points $x_1, x_2 \in \partial_- Y_3$, then these points are separated by a distance of at least $\delta$. It follows from this that there exists a point $y \in Y_3$ of infimal height with two distinct closest points $x_1, x_2 \in \partial_- Y_3$. For each $i$, let $\gamma_i : [0, h(y)] \to Y_3$ be a length minimising curve from $x_i$ to $y$.

We now show that $y$ does not lie on the lateral boundary $\partial_0 Y_3$ of $Y_3$. Indeed, otherwise, the closest point $x_3$ to $y$ in $\partial_- Y_2$ is an element of $Z := \partial\partial_- Y_2 = \partial\partial_- Y_3$. In particular, by minimality, $d(y, x_3) = d(y, x_1)$, so that $x_1$ is thus also an element of $Z$, for otherwise $\gamma_1$ would intersect $\partial_- Y_2$ at some point lying strictly closer to $y$ than $x_3$, which is absurd. Since the closest point to $y$ in $\partial_- Y_2$ is unique, it follows that $x_1 = x_3$, and a similar reasoning then shows that $x_2 = x_3 = x_1$, which is absurd. It follows that $y$ does not lie on the lateral boundary of $Y_3$, as asserted.

Since $y$ is trivially also not an element of the lower boundary of $Y_3$, it follows that $y$ is an interior point of this ideal convex cobordism. It then follows that $\dot{\gamma}_1(h(y)) = -\dot{\gamma}_2(h(y))$, for otherwise, by transversality, there would exist another point $y'$ of strictly lower height also having two distinct closest points in $\partial_- Y_3$, contradicting minimality of the height of $y$. The join of these two curves thus yields a geodesic segment $\gamma : [0, 2h(y)] \to Y_3$ starting and ending at two points of $\partial_- Y_3$. It remains only to show that this yields a contradiction. To this end, let $h_1$ and $h_2$ denote the respective height functions of $Y_1$ and $Y_2$. If $\gamma$ remains entirely in $Y_1$, then $(h_1 \circ \gamma)$ attains a maximum at some interior point, contradicting strict convexity. On the other hand, if $\gamma$ enters $Y_2$, then $(h_2 \circ \gamma)$ likewise attains a maximum at some interior point of $\gamma^{-1}(Y_2)$, also contradicting strict convexity. We have thus shown that every point $y$ of $Y_3$ has a unique closest point in $\partial_- Y_3$, as asserted, and this completes the proof. $\square$

**Lemma 4.4.2**

Let $X$ be a Cartan–Hadamard manifold of sectional curvature bounded above by $-1$. Let $(Y_1, \phi_1)$ and $(Y_2, \phi_2)$ be respectively a convex cobordism and an ideal convex cobordism in $X$. Suppose that the upper boundary of $(Y_1, \phi_1)$ is a reparametrisation of the lower boundary of $(Y_2, \phi_2)$ and let $(Y_3, \phi_3)$ denote their join. If $\mathcal{V}_2$ and $\mathcal{V}_3$ denote the respective leaf spaces of their vertical line foliations, then there exists a unique homeomorphism $\alpha : \mathcal{V}_3 \to \mathcal{V}_2$ such that, for all $L$ in $\mathcal{V}_3$, $\alpha(L)$ is asymptotic to $L$.

**Remark 4.4.1.** Observe that the lateral boundary of $\mathcal{V}_3$ naturally identifies with that of $\mathcal{V}_2$ and that $\alpha$ restricts to a homeomorphism over this set.

**Remark 4.4.2.** In particular, $\alpha$ also defines a natural homeomorphism from $\partial_\infty Y_3$ into $\partial_\infty Y_2$.

**Proof:** Indeed, let $\gamma : [0, \infty[ \to Y_3$ be a height-parametrised vertical line in $Y_3$. Let $t_0 > 0$ be such that $\gamma(t_0) \in \partial_+ Y_1 = \partial_- Y_2$. The restriction of $\gamma$ to $[t_0, \infty[$ is a complete, unbounded, unit speed geodesic ray in $Y_2$, and is thus asymptotic to a unique vertical line in $Y_2$, as desired. Conversely, let $\gamma : [0, \infty[ \to Y_2$ be a height parametrised vertical line in $Y_2$. for all positive $t$, let $x_t$ denote the closest point to $\gamma(t)$ in $\partial_- Y_3$. By compactness, $(x_t)_{t>0}$ has a limit point $x_\infty \in \partial_- Y_3$, and we verify that there exists a vertical line in $Y_3$ above $x_\infty$ asymptotic to $\gamma$. This yields the desired bijection $\alpha$. This function is trivially continuous, and since its domain is compact, it is a homeomorphism, as desired. $\square$

**4.5 - Existence and uniqueness.** We now prove the main result of this paper.

**Theorem 4.5.1, Existence and Uniqueness**

*Let $X$ be a Cartan–Hadamard manifold of sectional curvature pinched between $-K$ and $-1$ and let $\partial_\infty X$ denote its ideal boundary. Let $(Y,\phi)$ be an asymptotic Plateau problem in $\partial_\infty X$. If $(Y,\phi)$ is not of finite type 0, 1 or 2 then, for all $k \in ]0,1[$, there exists a unique quasicomplete ISC immersion $e : Y \to X$ of constant special lagrangian curvature equal to $k$ solving $(Y,\phi)$.*

The proof of Theorem 4.5.1 is given in the following series of lemmas. We first prove existence, which will require the following rather technical construction, illustrated in Figure 4.5.6. Choose $x_0 \in X$, and recall from Section 2.1 that this choice of base point defines a differential structure over $\partial_\infty X$. We furnish $Y$ with the unique differential structure making $\phi$ into a local diffeomorphism. Let $\Omega \subseteq Y$ be a relatively compact open subset with smooth boundary such that the restriction of $\phi$ to $\partial\Omega$ does not self-intersect tangentially. We first prove the existence of a solution to the asymptotic Plateau problem $(\Omega,\phi)$. Let $\Omega'$ be another relatively compact open subset of $Y$ with smooth boundary whose interior contains the closure of $\Omega$. For all $r$, let $S_r$ denote the geodesic sphere of radius $r$ about $x_0$ in $X$, let $\pi_r : \partial_\infty X \to S_r$ denote the projection along radial lines and denote $e_r := \pi_r \circ \phi$. For all $r$, let $(\hat{\Omega}_{r,\mathrm{icc}}, \hat{e}_{r,\mathrm{icc}})$ and $(\hat{\Omega}'_{r,\mathrm{icc}}, \hat{e}_{r,\mathrm{icc}})$ denote the respective ideal convex cobordisms over $(\Omega, \hat{e}_r)$ and $(\Omega', \hat{e}_r)$, and observe that the former naturally embeds in the latter. By comparison theory, for all $r$, $S_r$ has special lagrangian curvature bounded below at every point by $\coth(r)$. It follows by Theorem 4.3.1 that, for all $r$, there exists a unique convex cobordism $(\hat{\Omega}_{r,\mathrm{cc}}, \hat{e}_{r,\mathrm{cc}})$ whose upper boundary is a reparametrisation of the lower boundary of $(\hat{\Omega}_{r,\mathrm{icc}}, \hat{e}_{r,\mathrm{icc}})$ and whose lower boundary $(Z_r, f_r)$ has constant special lagrangian curvature equal to $k$. Let $(\hat{\Omega}_{r,\mathrm{sol}}, \hat{e}_{r,\mathrm{sol}})$ denote the ideal convex cobordism obtained by joining $(\hat{\Omega}_{r,\mathrm{cc}}, \hat{e}_{r,\mathrm{cc}})$ and $(\hat{\Omega}_{r,\mathrm{icc}}, \hat{e}_{r,\mathrm{icc}})$. Since $\hat{\Omega}_{r,\mathrm{icc}}$ naturally embeds inside $\hat{\Omega}'_{r,\mathrm{icc}}$ and $\hat{\Omega}_{r,\mathrm{sol}}$, we may join $\hat{\Omega}'_{r,\mathrm{icc}}$ and $\hat{\Omega}_{r,\mathrm{sol}}$ along this shared subset to obtain an ideal convex cobordism $(\hat{\Omega}'_{r,\mathrm{sol}}, \hat{e}'_{r,\mathrm{sol}})$, with piecewise smooth lower boundary, which contains both. Observe, finally, that $\partial_\infty \hat{\Omega}'_{r,\mathrm{icc}}$ naturally identifies with $\Omega'$ and that, by Lemma 4.4.2, $\partial_\infty \hat{\Omega}'_{r,\mathrm{sol}}$ also identifies with $\Omega'$.

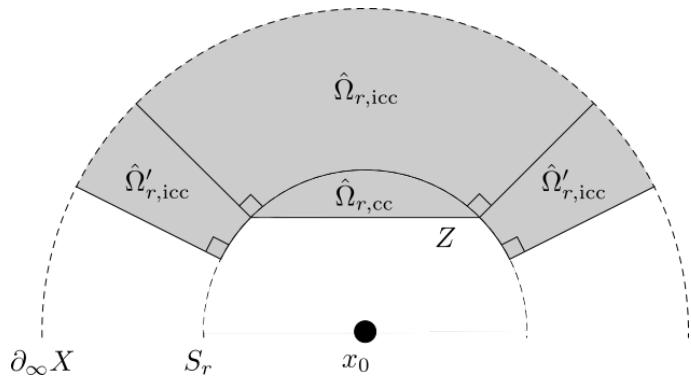

**Figure 4.5.6 - Constructing quasicomplete ISC immersions -** $\hat{\Omega}'_{r,\mathrm{sol}}$ is an ideal convex cobordism with piecewise smooth lower boundary obtained by joining $\hat{\Omega}_{r,\mathrm{sol}}$ to $\hat{\Omega}'_{r,\mathrm{icc}}$ along $\hat{\Omega}_{r,\mathrm{icc}}$.

**Lemma 4.5.2**

*There exists $r_0 > 0$ such that, for all $r > r_0$ and for every interior point $p$ of $\hat{\Omega}_{r,\mathrm{sol}}$, $T_{+,p}\hat{\Omega}_{r,\mathrm{sol}}$ is contained in the domain of the exponential map of $\hat{\Omega}'_{r,\mathrm{sol}}$.*

**Proof:** Indeed, choose $\epsilon > 0$. By Lemma 4.5.2, there exists $r_1 > 0$ such that, for every point $p$ of $\hat{\Omega}_{\epsilon,\mathrm{icc}}$ with $h(p) > r_1$, $T_{+,p}\hat{\Omega}_{\epsilon,\mathrm{icc}}$ is contained in the domain of the exponential map of $\hat{\Omega}'_{\epsilon,\mathrm{icc}}$. We now show that for $r > r_1$, $\hat{\Omega}_{r,\mathrm{sol}}$ has the desired property. Indeed, let $\gamma : [0,R[ \to \hat{\Omega}_{r,\mathrm{sol}}$ be a unit speed geodesic ray whose derivative is everywhere upward-pointing. If $\gamma$ remains in $\hat{\Omega}_{r,\mathrm{sol}}$, then it can be extended indefinitely. Otherwise, $\gamma$ meets the lateral boundary $\partial_0 \hat{\Omega}_{r,\mathrm{sol}}$ at some point $t$, say. Since the vertical line foliation of On the asymptotic Plateau problem in Cartan-Hadamard manifolds.$\partial_0\hat{\Omega}_{r,\text{sol}}$ concides with that of $\partial_0\hat{\Omega}_{r,\text{icc}}$, $\dot{\gamma}(t)$ is also upward-pointing in $\hat{\Omega}_{r,\text{icc}}$, and it follows by definition of $r_1$ that $\gamma$ extends indefinitely in $\hat{\Omega}'_{r,\text{icc}}$, as asserted. It follows that for every interior point $p$ of $\hat{\Omega}_{r,\text{sol}}$, $T_{+,p}\hat{\Omega}_{r,\text{sol}}$ is contained the domain of the exponential map of $\hat{\Omega}'_{r,\text{sol}}$, as desired. □

We henceforth suppose that $r > r_0$ and, for every interior point $p$ of $\hat{\Omega}_{r,\text{sol}}$, we denote by $(\text{AB}(p), \alpha(p))$ the asymptotic ball of $(\partial_\infty \hat{\Omega}'_{r,\text{sol}}, \partial_\infty \hat{e}'_{r,\text{sol}}) = (\Omega', \phi)$ centred on $p$.

Let $\pi'_{r,\text{icc}} : \Omega' \to \partial_-\hat{\Omega}'_{r,\text{icc}}$ and $\pi'_{r,\text{sol}} : \Omega' \to \partial_-\hat{\Omega}'_{r,\text{sol}}$ denote the respective projections along vertical vertical lines. These two projections coincide over the complement of $\Omega$ and $\pi'_{r,\text{icc}}$ is a diffeomorphism. Significantly, however, $\pi'_{r,\text{sol}}$ is not injective over the preimage of $\partial\partial_-\hat{\Omega}_{r,\text{sol}} = \partial Z_r$. We thus view the inverse $\psi_r$ of $\pi'_{r,\text{sol}}$ as a set-valued function having unique values away from $\partial Z_r$. Observe that, away from $\partial Z_r$, $\psi_r$ is a homeomorphism onto its image.

**Lemma 4.5.3**

*The restriction of $(\psi_r \circ \pi'_{r,\text{sol}})_{r>0}$ to $\partial\Omega$ converges uniformly to the identity with respect to the Hausdorff topology as $r$ tends to infinity.*

**Proof:** Suppose the contrary. By compactness of $\overline{\Omega}$, there exists a sequence $(r_m)_{m\in\mathbb{N}}$ converging to infinity, a sequence $(p_m)_{m\in\mathbb{N}}$ of points of $\partial\Omega$ converging to some point $p_\infty \in \partial\Omega$ and a sequence $(q_m)_{m\in\mathbb{N}}$ of points of $\overline{\Omega}$ converging to some point $q_\infty \neq p_\infty$ such that, for all $m$, $q_m$ is an element of $Q_m := (\psi_{r_m} \circ \pi'_{r_m,\text{sol}})(p_m)$. Let $(\text{AB}_\infty, \alpha_\infty)$ be an asymptotic ball of $(\Omega', \phi)$ about $p_\infty$. Observe that, for all $m$, $Q_m$ is connected, so that, upon reducing $\text{AB}_\infty$ if necessary, we may suppose that $q_\infty \in \partial\alpha_\infty(\text{AB}_\infty)$. For all $m$, let $\tilde{p}_m$ denote the point of $\hat{\Omega}_{r,\text{sol}}$ lying at height $(1/m)$ along the vertical line from $\pi'_{r_m,\text{sol}}(p_m)$ to $q_m$. For all $m$, let $(\text{AB}_m, \alpha_m)$ denote the asymptotic ball of $(\Omega', \phi)$ centred on $p_m$. Since $(\hat{e}'_{r_m,\text{sol}} \circ \tilde{p}_m)_{m\in\mathbb{N}}$ tends to $\phi(p_\infty)$ and $(\phi(q_m))_{m\in\mathbb{N}}$ tends to $\phi(q_\infty)$, it follows that the sequence $(\text{AB}^c_m)_{m\in\mathbb{N}}$ converges in the Hausdorff sense to the singleton $\{\phi(p_\infty)\}$. In particular, for sufficiently large $m$, the intersection of $\text{AB}_m$ with $\text{AB}_\infty$ is connected and their union covers the whole of $\partial_\infty X$. For all such $m$, the join of $\alpha_m$ with $\alpha_\infty$ thus defines a homeomorphism from $\partial_\infty X$ onto an open subset of $\Omega'$. $(Y, \phi)$ therefore contains an asymptotic Plateau problem of finite type 0. By Lemma 2.2.1, this is absurd, and the result follows. □

**Lemma 4.5.4**

*$(f_r \circ \pi'_{r,\text{sol}})_{r>0}$ subconverges in the $C^0_{\text{loc}}$ sense over $\Omega$ to a quasicomplete ISC immersion of constant special lagrangian curvature equal to $k$ solving the asymptotic Plateau problem $(\Omega, \phi)$.*

**Proof:** Indeed, for all $r$, $\psi_r$ restricts to a homeomorphism from the interior of $Z_r$ onto an open subset $\Omega_r$ of $\Omega$. In addition, for all $r$, the restriction of $(\phi \circ \psi_r)$ to $Z_r$ is none other than the asymptotic Gauss map of $f_r$. It follows that, for all $r$, $(f_r \circ \pi'_{r,\text{sol}})$ is a (non-quasicomplete) ISC immersion of constant special lagrangian curvature equal to $k$ solving $(\Omega_r, \phi)$. In addition, by Lemma 4.5.3, for every compact subset $K$ of $\Omega$, there exists $r_0 > 0$ such that, for all $r > r_0$, $K \subseteq \Omega_r$. Since $((f_r \circ \pi'_{r,\text{sol}})(\partial\Omega_r))_{r>0}$ tends towards a subset of $\partial_\infty X$, the result now follows by Theorem 4.2.1 and the subsequent remark. □

**Lemma 4.5.5**

*There exists a quasicomplete ISC immersion $e : Y \to X$ of constant special lagrangian curvature equal to $k$ solving the asymptotic Plateau problem $(Y, \phi)$.*

**Proof:** Let $(\Omega_m)_{m\in\mathbb{N}}$ be an exhausting, nested sequence of relatively compact open subsets of $Y$ with smooth boundary. By Lemma 4.5.4, for all $m$, there exists a quasicomplete ISC immersion $e_m : \Omega_m \to X$ of constant special lagrangian curvature equal to $k$ solving $(\Omega_m, \phi)$. By Theorem 4.2.1, the sequence $(e_m)_{m\in\mathbb{N}}$ subconverges to a quasicomplete ISC immersion $e : Y \to X$ of constant special lagrangian curvature equal to $k$ solving the asymptotic Plateau problem $(Y, \phi)$. This completes the proof. □

**Lemma 4.5.6**

*For all $k$, there exists at most one quasicomplete ISC immersion $e : Y \to X$ of constant special lagrangian curvature equal to $k$ solving the asymptotic Plateau problem $(Y, \phi)$.*

**Proof:** Suppose the contrary. Let $e_1, e_2 : Y \to M$ be two distinct solutions. Let $(\Omega_m)_{m\in\mathbb{N}}$ be an exhausting, nested sequence of relatively compact open subsets of $Y$ with smooth boundary. For all $m$, let $f_m : \Omega_m \to M$ be a quasicomplete ISC immersion of constant special lagrangian curvature equal to $k$ solving $(\Omega_m, \phi)$. By Lemma 3.4.2, for all $m$, $f_m$ is dominated by both $e_1$ and $e_2$. By Theorem 4.2.1, $(f_m)_{m \in \mathbb{N}}$ subconverges in the $C^0_{\text{loc}}$ sense towards a quasicomplete ISC immersion $f_\infty : Y \to M$ of constant special lagrangian curvature equal to $k$ which solves the asymptotic Plateau problem $(Y, \phi)$. We verify that $f_\infty$ is also dominated by both $e_1$ and $e_2$ so that, by Lemma 3.3.3, $e_1 = f_\infty = e_2$, as desired. $\square$

In order to prove Theorem 1.2.2, it remains only to exclude the exceptional cases.

**Theorem 4.5.7**

Let $X$ be a Cartan–Hadamard manifold of sectional curvature pinched between $-K$ and $-1$ and let $\partial_\infty X$ denote its ideal boundary. Let $(Y, \phi)$ be an asymptotic Plateau problem in $\partial_\infty X$. If $Y$ is of finite type 0, 1 or 2 then, for all $k \in\, ]0, 1[$, there exists no quasicomplete ISC immersion $e : Y \to X$ of constant special lagrangian curvature equal to $k$ solving $(Y, \phi)$.

**Proof:** The cases where $(Y, \phi)$ is of finite type 0 or 1 and $X$ is 3-dimensional are addressed by Theorem 7.4.1 of [25], and the same argument also applies to the higher dimensional case. The case where $(Y, \phi)$ is of finite type 2 and $\phi$ is a homeomorphism is also addressed by this result. It remains only to address the case where $(Y, \phi)$ is of finite type 2 and $\phi$ is a non-trivial cover, which is only possible when $X$ is 3-dimensional. Suppose that $\phi$ is an $m$-fold cover, for some $m \in \mathbb{N} \cup \{\infty\}$. Let $\text{Hor} : SX \to \partial_\infty X$ denote the horizon map of $X$, let $p, q \in \partial_\infty X$ denote the two points of the complement of $\phi(Y)$, let $\Gamma$ denote the unique geodesic in $X$ joining $p$ and $q$, and let $S\Gamma$ denote the circle bundle of unit normal vectors over $\Gamma$. There trivially exists a unique $m$-fold cover $\hat{e} : Y \to S\Gamma$ such that $\text{Hor} \circ \hat{e} = \phi$. Bearing in mind Theorem 4.1.2, $(Y, \hat{e})$ may be viewed as a degenerate solution of the asymptotic Plateau problem $(Y, \phi)$. We verify that Lemma 4.5.6 also applies in this case, from which it follows that the asymptotic Plateau problem $(Y, \phi)$ has no non-degenerate solutions, and this completes the proof. $\square$

**4.6 - Dynamical stability.** We now prove stability of $k$-hypersurface laminations.

**Theorem 1.1.1**

Let $X$ and $X'$ be Cartan–Hadamard manifolds of sectional curvature pinched between $-K$ and $-1$ for some $K \geqslant 1$. For all $k, k' \in\, ]0, 1[$ and for every homeomorphism $\psi : \partial_\infty X \to \partial_\infty X'$, there exists a unique homeomorphism $\Psi : \overline{\mathcal{S}}_k(X) \to \overline{\mathcal{S}}_{k'}(X')$ such that the following diagram commutes.

$$\begin{array}{ccc} \overline{\mathcal{S}}_k(X) & \xrightarrow{\psi} & \overline{\mathcal{S}}_{k'}(X') \\ \Phi \downarrow & & \downarrow \Phi \\ \partial_\infty X & \xrightarrow{\Psi} & \partial_\infty X'. \end{array}$$

Furthermore, $\Psi$ maps $\partial \mathcal{S}(X)$ into $\partial S(X')$.

We first prove two technical lemmas. Let $X$ be a Cartan–Hadamard manifold of sectional curvature pinched between $-K$ and $-1$ for some $K \geqslant 1$, and let $\partial_\infty X$ denote its ideal boundary. Choose $k \in\, ]0, 1[$, and let $(Y_m, e_m, y_m)_{m \in \mathbb{N}}$ be a sequence of pointed $k$-hypersurfaces in $X$ converging towards $(Y_\infty, e_\infty, y_\infty)$. For all $m \in \mathbb{N} \cup \{\infty\}$, let $\hat{e}_m$ denote the Gauss lift of $e_m$, and let $\phi_m$ denote its asymptotic Gauss map.

**Lemma 4.6.1**

Let $(Y\phi, , y)$ be a pointed asymptotic Plateau problem in $\partial_\infty X$ that is not of finite type 0, 1 or 2. If, for all $m$, there exists an injective homeomorphism $\alpha_m : Y \to Y_m$ such that $\alpha_m(y) = y_m$ and $\phi = \phi_m \circ \alpha_m$, then there exists an injective homeomorphism $\alpha_\infty : Y \to Y_\infty$ such that $\alpha_\infty(y) = y_\infty$ and $\phi = \phi_\infty \circ \alpha_\infty$.

**Proof:** Solving the asymptotic Plateau problem $(Y, \phi)$ yields a natural differential structure over $Y$. Let $\Omega$, $\Omega'$ be relatively compact open subsets of $Y$ with smooth boundaries such that

$$y \in \Omega \subseteq \overline{\Omega} \subseteq \Omega'.$$

Since $\phi^{-1}(\{\phi(y)\})$ is discrete, we may suppose, in addition, that $\phi(y) \notin \phi(\partial\Omega)$.

Let $f : \Omega \to X$ be the $k$-surface solving the asymptotic Plateau problem $(\Omega, \phi)$. For all $m$, let $(\hat{Y}_m, \hat{e}_m)$ denote the end of $(Y_m, e_m)$, as in Lemma and Definition 2.8.1. By Lemmas 2.4.6 and 3.4.2, for all $m$, $(Y_m, e_m)$ dominates $(f, \Omega)$ so that there exists a smooth embedding $\hat{f}_m : \Omega \to \hat{Y}_m$ satisfying $f = \hat{e}_m \circ \hat{f}_m$.

We now claim that there exists $\epsilon > 0$ such that, for all $m$, $\hat{f}_m(\Omega)$ lies above the level of $\hat{Y}_m$ at height $\epsilon$. Indeed, let $f' : \Omega' \to X$ be the $k$-surface solving the asymptotic Plateau problem $(\Omega', \phi)$ and let $(\hat{\Omega}', \hat{f}')$ denote its end. As above, there exists a smooth embedding $\hat{f} : \Omega \to \hat{\Omega}'$ such that $f = \hat{f}' \circ \hat{f}$. Since $\Omega$ is a relatively compact subset of $\Omega'$, there exists $\epsilon > 0$ such that $\hat{f}(\Omega)$ lies above the level of $\hat{\Omega}'$ at height $\epsilon$. Since, for all $m$, $(Y_m, e_m)$ also dominates $(\Omega', f')$, it follows that this value of $\epsilon$ has the desired properties.

For all $m$, let $z_m$ denote the point of $\Omega$ whose image under $\hat{f}_m$ lies above $y_m$. We claim that $(z_m)_{m \in \mathbb{N}}$ lies within a compact subset of $\Omega$. Indeed, otherwise, we may suppose that $(z_m)_{m \in \mathbb{N}}$ converges to some point $z_\infty$ of $\partial\Omega$. By Lemma 3.4.1, $(f(z_m))_{m \in \mathbb{N}}$ converges to $\phi(z_\infty)$. However, for all $m$, $f(z_m)$ lies on the geodesic ray joining $e_m(y_m)$ and $\phi(y)$, which is absurd, since $(e_m(y_m))_{m \in \mathbb{N}}$ is bounded and $\phi(z_\infty) \neq \phi(y) \notin \phi(\partial\Omega)$.

We now claim that, for every compact subset $K$ of $\Omega$, there exists $R$ such that, for all $m$, $\hat{f}_m(K)$ lies above the ball of radius $R$ about $y_m$ in $(Y_m, e_m)$. Indeed, by the preceeding argument, there exists $R' > 0$ such that, for all $m$, $K$ is contained in the ball of radius $R'$ about $z_m$ in $(\Omega, f)$. Since projection along vertical lines onto every level is distance decreasing, it follows by (2.38) that $R := \sqrt{2}\cosh_K(\epsilon)R'$ has the desired properties.

Upon taking limits, we now see that $f(K)$ is a graph over some subset of $B_R(y_\infty)$. Upon taking an exhaustion of $\Omega$ by compact sets, we then see that $(\Omega, f)$ is dominated by $(Y_\infty, e_\infty)$. There therefore exists an injective homeomorphism $\alpha_\infty : \Omega \to Y_\infty$ such that $\alpha_\infty(y) = y_\infty$ and $\phi = \phi_\infty \circ \alpha_\infty$. The result now follows upon covering $Y$ by a nested union of such open sets. $\square$

Now let $X'$ be another Cartan–Hadamard manifold also of sectional curvature pinched between $-K$ and $-1$, let $\partial_\infty X'$ denote its ideal boundary, and let $\alpha : \partial_\infty X \to \partial_\infty X'$ be a homeomorphism.

**Lemma 4.6.2**

Let $(Y_m, \phi_m, y_m)_{m \in \mathbb{N}}$ be a sequence of pointed asymptotic Plateau problems in $\partial_\infty X$. For all $m$, let $(Y_m, e_m)$ and $(Y_m, f_m)$ solve $(Y_m, \phi_m)$ and $(Y_m, \alpha \circ \phi_m)$ respectively. Then $(e_m(y_m))_{m \in \mathbb{N}}$ is bounded in $X$ if and only if $(f_m(y_m))_{m \in \mathbb{N}}$ is bounded in $Y$.

**Proof:** Suppose the contrary, so that, without loss of generality, $(e_m(y_m))_{m \in \mathbb{N}}$ is bounded in $X$ but $(f_m(y_m))_{m \in \mathbb{N}}$ is unbounded in $X'$. We may suppose that $(\phi_m(y_m))_{m \in \mathbb{N}}$ converges to $y_\infty$, say, in $\partial_\infty X$, that $(e_m(y_m))_{m \in \mathbb{N}}$ converges to $a_\infty$, say, in $X$, and that $(f_m(y_m))_{m \in \mathbb{N}}$ converges to $b_\infty$, say, in $\partial_\infty X'$. We may also suppose that, for all $m$, $(Y_m, \phi_m)$ contains a fixed asymptotic ball about $y_m$ so that, by Lemma 4.2.2, $b_\infty$ is distinct from $\alpha(y_\infty)$.

Let $\gamma : \mathbb{R} \to X'$ denote the geodesic joining $b_\infty$ to $\alpha(y_\infty)$. Observe that, for all $t$, and for all sufficiently large $m$, the asymptotic ball $\text{AB}(\dot{\gamma}(t))$ is contained within $(Y_m, \alpha \circ \phi_m)$, so that $\alpha^{-1}(\text{AB}(\dot{\gamma}(t)))$ is contained within $(\phi_m, Y_m)$. However, by Theorems 4.1.1 and 4.1.2, we may suppose that $(Y_m, e_m, y_m)_{m \in \mathbb{N}}$ converges to a $k$-hypersurface $(Y_\infty, e_\infty, y_\infty)$ with asymptotic Gauss map $\psi_\infty$. By Lemma 4.6.1, $(Y_\infty, \psi_\infty)$ contains $\alpha^{-1}(\text{AB}(\dot{\gamma}(t)))$ for all $t$. Upon letting $t$ tend to $-\infty$, it follows that $(\psi_\infty, Y_\infty)$ contains $\partial_\infty X \setminus \{\alpha^{-1}(b_\infty)\}$ so that, by Lemma 2.2.1, $(\psi_\infty, Y_\infty)$ is of finite type 1, which is absurd by Theorem 1.2.2. The result follows. $\square$

**Proof of Theorem 1.1.1:** The construction of $\Psi$ is elementary. Indeed, choose $(Y, e, y) \in \mathcal{S}_k(X)$. Let $\phi_e$ denote its asymptotic Gauss map. By Theorem 1.2.2, $(Y, \phi_e)$ is not of finite type 0, 1 or 2, and therefore neither is $(Y, \psi \circ \phi_e)$. It follows by Theorem 1.2.2 again that there exists a unique $k'$-hypersurface $e' : Y \to X'$ solving the asymptotic Plateau problem $(Y, \psi \circ \phi_e)$. We then define

$$\Psi(Y, e, y) := (Y, e, y).$$

The definition of $\Psi$ over $\partial\mathcal{S}(X)$ is straightforward. By construction, $\Psi$ is a bijection with the required properties. It remains only to show that $\Psi$ is a homeomorphism.

By symmetry, it suffices to show that $\Psi$ is continuous. Thus, let $(Y_m, \hat{e}_m, y_m)_{m \in \mathbb{N}}$ be a sequence in $\overline{\mathcal{S}}_k(X)$ converging to $(Y_\infty, \hat{e}_\infty, y_\infty)$. For all $m \in \mathbb{N} \cup \{\infty\}$, let $\hat{f}_m : Y_m \to UX'$ be such that

$$\phi_m := \phi_{\hat{f}_m} = \psi \circ \phi_{\hat{e}_m}.$$

By Lemma 4.6.2, the sequence $(\hat{f}_m(y_m))_{m\in\mathbb{N}}$ is bounded in $X'$. It follows by Theorems 4.1.1 and 4.1.2 that $(Y_m, \hat{f}_m, y_m)_{m\in\mathbb{N}}$ is relatively compact in $\overline{\mathcal{S}}_k(X')$. Let $(Y'_\infty, \hat{f}'_\infty, y'_\infty)$ be an accumulation point of this set and denote

$$\phi'_\infty := \phi_{\hat{f}'_\infty}.$$

We claim that there exists an injective local homeomorphism $\alpha : Y'_\infty \to Y_\infty$ such that

$$\phi'_\infty = \phi_\infty \circ \alpha.$$

Indeed, let $(\Phi_m)_{m\in\mathbb{N}}$ be a sequence of convergence maps of $(Y_m, \hat{f}_m, y_m)_{m\in\mathbb{N}}$ with respect to $(Y'_\infty, \hat{f}'_\infty, y'_\infty)$. Let $\Omega$ and $\Omega'$ be relatively compact neighbourhoods of $y'_\infty$ in $Y'_\infty$ such that $\overline{\Omega} \subseteq \Omega'$. By definition $(\phi_m \circ \Phi_m)_{m\in\mathbb{N}}$ converges to $\phi'_\infty$ locally uniformly over $\Omega'$. Upon increasing $M$ if necessary, we may therefore suppose that, for all $m \geqslant M$, there exists an injective local homeomorphism $\alpha_m : \Omega \to Y_m$ such that $\alpha_m(y'_\infty) = y_m$ and $\phi'_\infty = \phi_m \circ \alpha_m$. By Lemma 4.6.1, there exists an injective local homeomorphism $\alpha : \Omega \to Y_\infty$ such that $\alpha(y'_\infty) = y_m$ and $\phi'_\infty = \phi_m \circ \alpha$. Taking an exhaustion of $Y'_\infty$ by such open sets yields an injective homeomorphism $\alpha : Y'_\infty \to Y_\infty$ with the desired properties, proving the assertion. By symmetry, there likewise exists an injective local homeomorphism $\beta : Y_\infty \to Y'_\infty$ such that

$$\phi_\infty = \phi'_\infty \circ \beta.$$

It follows that $\alpha$ and $\beta$ are homeomorphisms. In particular, the sequence $(Y_m, \hat{f}_m, y_m)_{m\in\mathbb{N}}$ has $(Y_\infty, \hat{f}_\infty, y_\infty)$ as its unique accumulation point, and thus converges towards this point. Continuity follows, and this completes the proof. $\square$

## A - Comparison theory.

Comparison theory makes it useful to use a consistent terminology for trigonometric and hyperbolic functions. Thus, for $k \in \mathbb{R}$, we define the functions $\sinh_k, \cosh_k : \mathbb{R} \to \mathbb{R}$ to be the unique solutions of the ordinary differential equation

$$\phi'' - k\phi = 0, \tag{A.1}$$

with initial conditions

$$\sinh_k(0) = \cosh'_k(0) = 0 \text{ and } \sinh'_k(0) = \cosh_k(0) = 1. \tag{A.2}$$

We then define, for all $k$ and for all $x$,

$$\begin{aligned}\tanh_k(x) &:= \frac{\sinh_k(x)}{\cosh_k(x)} \text{ and} \\ \coth_k(x) &:= \frac{\cosh_k(x)}{\sinh_k(x)}.\end{aligned} \tag{A.3}$$

These are the only such functions that we will require in this paper.

## B - Omori's maximum principle.

In this appendix, we sketch a proof of Omori's maximum principle in the form used in the proof of Lemma 3.3.2. Our approach is slightly more direct than that given by Omori in [28]. We will require the following useful formula. Let $X$ be a riemannian manifold, and let $Y$ be an embedded hypersurface in $X$. Let $\nu : Y \to TX$ be a unit normal vector field over $Y$ and let $II$ denote its corresponding second fundamental form. Let $\nabla$, Hess and $\overline{\nabla}$, $\overline{\text{Hess}}$ denote the respective gradient and Hessian operators of $Y$ and $X$. Then, for all twice differentiable $f : X \to \mathbb{R}$,

$$\text{Hess}(f) = \overline{\text{Hess}}(f)|_{TY} - \langle \overline{\nabla} f, \nu \rangle II. \tag{B.1}$$

**Theorem B.1**

*Let $X$ be a riemannian manifold with sectional curvature bounded below. Let $f : X \to \mathbb{R}$ be a $C^2$ function. Suppose that*

$$C := \sup_{x \in X} f(x) < \infty, \tag{B.2}$$

*and that there exists $\epsilon > 0$ such that, for all $x \in X$ with $f(x) > C - \epsilon$, the closed ball $\overline{B}_\epsilon(x)$ is complete. Then, for all $\delta > 0$, there exists $x \in X$ such that*

$$\begin{aligned} f(x) &> C - \delta, \\ \|Df(x)\| &< \delta, \text{ and} \\ \mathrm{Hess}(f)(x) &< \delta \mathrm{Id}. \end{aligned} \tag{B.3}$$

**Proof:** Without loss of generality, we may suppose that $C = 0$ and that $\epsilon = 1$. Let $\Gamma$ denote the graph of $f$ in $X \times \mathbb{R}$. Let $\pi : \Gamma \to X$ denote the projection onto the first factor, let $\mathrm{I}^\Gamma$ denote the first fundamental form of $\Gamma$, let $\nu^\Gamma$ denote its downward pointing unit normal vector field, and let $\mathrm{II}^\Gamma$ denote its corresponding second fundamental form. The result will follow upon determining estimates for $\mathrm{II}^\Gamma$ at suitably chosen points.

First observe that $\Gamma$ is the level set of the function

$$\hat{f}(x, t) := f(x) - t.$$

From this, it readily follows that, for all $(x, t) \in \Gamma$,

$$\nu^\Gamma(x, t) = (1 + \|\nabla f(x)\|^2)^{-1/2}(\nabla f(x), -1),$$

and, by (B.1),

$$\mathrm{II}^\Gamma(x, t) = (1 + \|\nabla f(x)\|^2)^{-1/2} \pi^* \mathrm{Hess}(f)(x).$$

Now choose $0 < \delta < 1$, choose $x_0 \in X$ such that $f(x_0) > -\delta$ and choose $a > 0$ such that

$$a < \frac{1}{2}\left(\frac{1}{\delta} - \delta\right). \tag{B.4}$$

Let $\hat{d}$ denote the distance in $X \times \mathbb{R}$ to $(x_0, a)$ and let $d$ denote the distance in $X$ to $x_0$. Let $(y, b)$ be a point in $\Gamma$ minimising $\hat{d}$. It follows from (B.4) that $y \in B_1(x_0)$. Let $\hat{\gamma} : [0, l] \to X \times \mathbb{R}$ be a unit speed parametrised length-minimising geodesic from $(x_0, a)$ to $(y, b)$. In Lemma 5 of [28], Omori shows that

$$\hat{\gamma}(t) = (\gamma(t), tb/l + (l - t)a/l),$$

for some geodesic $\gamma : [0, l] \to X$. In addition, $y$ is not conjugate to $x_0$ along $\gamma$ and, for all $(y', b')$ in a neighbourhood of $(y, b)$,

$$\hat{d}(y', b')^2 = d(y')^2 + (a - b')^2.$$

Since $(y, b)$ minimises $\hat{d}$ over $\Gamma$, the sphere of radius $\hat{d}(y, b)$ about $(x_0, a)$ is an interior tangent to $\Gamma$ at this point. It follows that

$$\nu^\Gamma(y, b) = (\nabla \hat{d})(y, b) = \frac{1}{\hat{d}(y, b)}(d(y)(\nabla d)(y), -(a - b)).$$

Comparing the last components of these vectors yields

$$1 + \|(\nabla f)(y)\|^2 \leqslant \left(1 + \frac{1}{a}\right)^2.$$

In particular,
$$\|(\nabla f)(y)\|^2 \leqslant \frac{2}{a} + \frac{1}{a^2}.$$

Likewise,
$$\begin{aligned} 2a\mathrm{II}^\Gamma(y,b) &\leqslant 2\hat{d}(y,b)\mathrm{II}^\Gamma(y,b) \\ &= \langle (\nabla \hat{d}^2)(y,b), \nu^\Gamma(y,b)\rangle \mathrm{II}^\Gamma(y,b) \\ &\leqslant \mathrm{Hess}(\hat{d}^2)(y,b)|_{T\Gamma}, \end{aligned}$$

where the last inequality follows from (B.1) again.

Finally, by comparison theory, there exists $B > 0$, which only depends on the lower bound for the sectional curvature of $X$, such that
$$\mathrm{Hess}(d^2)(y) \leqslant B\mathrm{Id},$$

so that
$$\mathrm{Hess}(\hat{d}^2)(y,b) = \mathrm{Hess}(d^2)(y) + 2dt^2 \leqslant B\mathrm{Id}.$$

Hence
$$\begin{aligned} \mathrm{Hess}(f)(y) &= (1 + \|\nabla f(y)\|^2)^{1/2} \pi_* \mathrm{II}^\Gamma(y,b) \\ &\leqslant \left(1 + \frac{1}{a}\right) \frac{B}{2a} \pi_* \mathrm{I}^\Gamma(y,b) \\ &\leqslant \left(1 + \frac{1}{a}\right)^3 \frac{B}{2a} \mathrm{Id}, \end{aligned}$$

and the result follows upon choosing $\delta$ sufficiently small and $a$ sufficiently large. $\square$

## C - Special legendrian submanifolds and the holographic principle.

When expressed in terms of special legendrian submanifolds of the unit sphere bundle, Theorem 1.2.2 may be viewed as another manifestation of the holographic principle. This will be the content of the result discussed in this section. Let $X$ be a $(d+1)$-dimensional Cartan–Hadamard manifold of sectional curvature pinched between $-K$ and $-1$ for some $K \geqslant 1$, let $\partial_\infty X$ denote its ideal boundary, let $SX$ denote its unit sphere bundle and let $\pi : SX \to X$ denote the canonical projection. The total space of $SX$ carries a natural 1-parameter family of "almost Calabi-Yau" structures, as we now proceed to show.

First, upon identifying $T^*X$ with $TX$ via the riemannian metric, the canonical Liouville form of $T^*X$ restricts to a contact form over $SX$. The contact subbundle of $TSX$, that is, the kernel of this form, also has the following more explicit description. First, the Levi-Civita connection of $X$ yields a canonical splitting

$$TSX = HSX \oplus VSX, \tag{C.1}$$

of $TSX$ into a horizontal and vertical component. The subbundle $HSX$ is naturally isomorphic to $\pi^*TX$, whilst the subbundle $VSX$ is naturally isomorphic to the subbundle of $\pi^*TX$ whose fibre at any point $\xi$ is the orthogonal complement of $\xi$ in $T_{\pi(\xi)}X$. With these identifications, we define the codimension 1 subbundle $WSX \subseteq TSX$ such that, for all $\xi \in SX$,

$$W_\xi SX := \langle \xi \rangle^\perp \oplus \langle \xi \rangle^\perp, \tag{C.2}$$

and we verify that $WSX$ is none other than the above defined contact subbundle of $SX$.

The above explicit description yields a number of natural geometric structures over the contact subbundle. First, we define the **Minkowski metric** $m$ over this subbundle such that, with respect to the identification (C.2),

$$m((\xi,\nu),(\xi,\nu)) := 2\langle \xi, \nu \rangle. \tag{C.3}$$

Next, for $\lambda \in \mathbb{R}$, we define the $\lambda$-**twisted Sasaki metric** $g_\lambda$ and the $\lambda$-**twisted almost complex structure** $J_\lambda$ over this subbundle such that, with respect to this identification,

$$g_\lambda((\xi,\nu),(\xi,\nu)) := \lambda\|\xi\|^2 + \lambda^{-1}\|\nu\|^2, \text{ and} \tag{C.4}$$
$$J_\lambda(\xi,\nu) := (-\lambda^{-1}\nu, \lambda\xi).$$

For all $\lambda$, $J_\lambda$ is trivially compatible with $g_\lambda$. For any lagrangian subspace $L$ of $WSX$, we define $\mathrm{dVol}_{\lambda,L}$ to be the volume form of the restriction of $g_\lambda$ to this subspace. For all $\lambda$, we then define the almost complex $d$-form $\Omega_\lambda$ to be the unique $J_\lambda$-complex $d$-form over $WSX$ whose restriction to every fibre $H$ of $HSX$ coincides with $\mathrm{dVol}_{\lambda,H}$. For all $\lambda$, for all $x \in SX$, and for any lagrangian subspace $L$ of $W_\xi SX$, there exists a unique angle $\theta$ such that

$$\Omega_\lambda|_L = e^{i\theta}\mathrm{dVol}_{\lambda,L}. \tag{C.5}$$

We call $\theta$ the $\lambda$-**special lagrangian angle** of $L$ (c.f. [18]).

We say that an immersed submanifold $(Y,\hat{e})$ of $SX$ is **positive** whenever the restriction of $m$ to this submanifold is everywhere positive definite. We say that $(Y,\hat{e})$ is $(\lambda,\theta)$-**special legendrian** whenever it is legendrian with constant $\lambda$-special lagrangian angle equal to $\theta$. Now let $(Y,e)$ be an immersed submanifold in $X$ with unit normal vector field $\nu_e$. In [35] we show that $(Y,e)$ is infinitesimally strictly convex if and only if $\hat{e} := \nu_e$ is positive. Likewise, for all $\lambda$, we show that this hypersurface has constant special lagrangian curvature equal to $\lambda/\tanh((d-1)\pi/2d)$ if and only if $\hat{e} := \nu_e$ is $(\lambda,(d-1)\pi/2)$-special legendrian. Theorem 1.2.2 yields the following result concerning the prescription of complete special legendrian submanifolds of $SX$.

**Theorem C.1**

*Let $X$ be a $(d+1)$-dimensional Cartan–Hadamard manifold of sectional curvature pinched between $-K$ and $-1$ for some $K \geqslant 1$, let $\partial_\infty X$ denote its ideal boundary and let $SX$ denote its unit sphere bundle. For every asymptotic Plateau problem $(Y,e)$ in $\partial_\infty X$ not of finite type 0 or 1, and for all $0 < \lambda < \tanh((d-1)\pi/2d)$, there exists a unique, complete, positive $(\lambda,(d-1)\pi/2)$-special legendrian immersion $\hat{e} : Y \to SX$ such that*

$$\mathrm{Hor} \circ \hat{e} = \phi. \tag{C.6}$$

**Remark C.1.** In the case where $(Y,\phi)$ is of finite type 2, the desired special legendrian immersion is given by a tube over the unique geodesic joining the two points in the complement of $\phi(Y)$.

## D - The ideal boundary functor and its right inverse.

Our results also yield an intriguing generalisation to the non-constant curvature case of the construction [20] of Kulkarni–Pinkall which we now describe. In what follows, we will require the concepts of developed Cartan–Hadamard ends and their ideal boundaries, introduced in the second half of Chapter 2.

Let $X$ be a Cartan–Hadamard manifold with ideal boundary $\partial_\infty X$. Recall that the ideal boundary operator $\partial_\infty$ defines a covariant functor from the category of developed Cartan–Hadamard ends in $X$ into the category of asymptotic Plateau problems in $\partial_\infty X$. In [37], we show that, in the case where $X = \mathbb{H}^{d+1}$ is $(d+1)$-dimensional hyperbolic space, the construction [20] of Kulkarni–Pinkall may be understood as a right inverse $\mathcal{H}_{-1}$ of this functor. Indeed, given an asymptotic Plateau problem $(Y,\phi)$ in $\partial_\infty \mathbb{H}^{d+1}$ not of finite type 0 or 1, $(\mathcal{H}_{-1}Y, \mathcal{H}_{-1}\phi)$ is simply the unique maximal developed hyperbolic end with ideal boundary $(Y,\phi)$.

Theorem 1.2.2 allows us to extend Kulkarni–Pinkall's construction as follows. Suppose that $X$ has sectional curvature pinched between $-K$ and $-1$, for some $K \geqslant 1$. Let $(Y,\phi)$ be an asymptotic Plateau problem in $\partial_\infty X$ not of finite type 0, 1 or 2. For $0 < k < 1$, let $e_k : Y \to X$ denote the unique quasicomplete ISC immersion of constant special lagrangian curvature equal to $k$ with asymptotic Gauss map $\phi$, and let $(\hat{Y}_k, \hat{e}_k)$ denote its end, as defined in Section 2.8. An elementary modification of the uniqueness argument of Lemma 4.5.6 shows that, for all $k < k'$, there exists a canonical isometric embedding $\alpha : \hat{Y}_{k'} \to \hat{Y}_k$ such that $\hat{e}_{k'} = \hat{e}_k \circ \alpha$. The direct limit of the family $(\hat{Y}_k, \hat{e}_k)_{k \in ]0,1[}$ then yields a developed Cartan–Hadamard end, which we denote by $(\mathcal{H}Y, \mathcal{H}\phi)$ and whose ideal boundary is canonically isomorphic to $(Y,\phi)$.

When $(Y,\phi)$ is of finite type 2, we may also define $(\mathcal{H}Y, \mathcal{H}\phi)$ by taking a suitable cover of the complement of the unique geodesic in $X$ joining the two points in the complement of $\phi(Y)$. The reader may then verify that $\mathcal{H}$ is functorial and that, when $X$ is hyperbolic, $\mathcal{H}$ coincides with $\mathcal{H}_{-1}$. We have thus proven the following result.

**Theorem D.1**

Let $X$ be a Cartan–Hadamard manifold of sectional curvature pinched between $-K$ and $-1$ for some $K \geqslant 1$. There exists a functor $\mathcal{H}$ sending the category of asymptotic Plateau problems in $\partial_\infty X$ which are not of finite type $0$ or $1$ into the category of developed Cartan–Hadamard ends in $X$. Furthermore, this functor is a right inverse of $\partial_\infty$ and, when $X = \mathbb{H}^{d+1}$ is $(d+1)$-dimensional hyperbolic space, it coincides with $\mathcal{H}_{-1}$.

**Remark D.2.** Consider two asymptotic Plateau problems $(Y, \phi)$ and $(Y', \phi')$ in $\partial_\infty X$, not of finite type $0$ or $1$, such that $(Y', \phi')$ is contained in $(Y, \phi)$. When $X$ is hyperbolic, by maximality, every Cartan–Hadamard end in $X$ with ideal boundary $(Y', \phi')$ naturally embeds in $\mathcal{H}_{-1}Y$. It would be interesting to know whether this maximality continues to hold in the general case. By construction, for all $0 < k < 1$, the unique quasicomplete ISC immersion $e : Y' \to X$ of constant special lagrangian curvature equal to $k$ solving $(Y', \phi')$ naturally lifts to an embedding of $Y'$ into $\mathcal{H}Y$. This constitutes a partial answer to the question of maximality.

## E - Bibliography.

[1] Alvarez S., Lowe B., Smith G., Foliated Plateau problems and asymptotic counting of surface subgroups, arXiv:2212.13604

[2] Anderson M. T., The Dirichlet problem at infinity for manifolds of negative curvature, *J. Diff. Geom.*, **18**, (1983), 701–721

[3] Ballmann W., Gromov M., Schroeder V., *Manifolds of Nonpositive Curvature*, Birkhäuser Verlag, (1985)

[4] Barbot T., Béguin F., Zeghib A., Prescribing Gauss curvature of surfaces in 3-dimensional spacetimes, Application to the Minkowski problem in the Minkowski space, *Ann. Inst. Fourier*, **61**, no. 2, (2011), 511–591

[5] Bonsante F., Mondello G., Schlenker J. M., A cyclic extension of the earthquake flow I., *Geom. Topol.*, **17**, no. 1, (2013), 157–234

[6] Bonsante F., Mondello G., Schlenker J. M., A cyclic extension of the earthquake flow II., *Ann. Sci. Éc. Norm. Supér.*, **48**, no. 4, (2015), 811–859

[7] Caffarelli L., Kohn J. J., Nirenberg L., Spruck J., The Dirichlet problem for nonlinear second-order elliptic equations. II. Complex MongeAmpère, and uniformly elliptic, equations, *Comm. Pure Appl. Math.*, **38**, no. 2, (1985), 209–252

[8] Calabi E., An extension of E. Hopf's maximum principle with an application to riemannian geometry, *Duke Math. J.*, **25**, no. 1, (1958), 45–56

[9] Calegari D., Marques F. C., Neves A., Counting minimal surfaces in negatively curved 3-manifolds, to appear in *Duke Math. Journ.*

[10] Marques F. C., Neves A., Min-max theory and the Willmore conjecture, *Ann. Math.*, **179**, 683–782

[11] Marques F. C., Neves A., Song A., Equidistribution of minimal hypersurfaces for generic metrics, *Invent. Math.*, **216**, 423–443

[12] Crandall M. G., Ishii H., Lions P. L., User's guide to viscosity solutions of second order partial differential equations, *Bull. Amer. Math. Soc.*, **27**, (1992), 1–67

[13] Gromov M., Sign and geometric meaning of curvature, *Rendiconti del Seminario Matematico e Fisico di Milano*, **61**, (1991), 9–123

[14] Gromov M., Foliated Plateau problem, part I: minimal varieties, *G.A.F.A.*, **1**, no. 1, (1991), 14–79

[15] Gromov M., Foliated Plateau problem, part II: harmonic maps of foliations, *G.A.F.A.*, **1**, no. 1, (1991), 253–320

[16] Gromov M., Three remarks on geodesic dynamics and fundamental group, *Enseign. Math.*, **46**, 391402, (2000)

[17] Guan B., Spruck J., Hypersurfaces of Constant Mean Curvature in Hyperbolic Space with Prescribed Asymptotic Boundary at Infinity, *Amer. J. Math.*, **122**, no. 5, 1039–1060, (2000)

[18] Harvey F. R., Lawson H. B., Calibrated Geometries, *Acta Math.*, **148**, (1982), 47–157

[19] Harvey F. R., Lawson H. B., Pseudoconvexity for the special Lagrangian potential equation, *Calc. Var. PDEs.*, **60**, (2021)

[20] Kulkarni R.S., Pinkall U., A canonical metric for Möbius structures and its applications, *Math. Z.*, **216**, no. 1, (1994), 89–129

[21] Labourie F, Immersions isométriques elliptiques et courbes pseudoholomorphes, *J. Diff. Geom.*, **30**, (1989), 395-424

[22] Labourie F., Problème de Minkowski et surfaces à courbure constante dans les variétés hyperboliques, *Bull. Soc. Math. Fr.*, **119**, (1991), 307–325

[23] Labourie F., Métriques prescrites sur le bord des variétés hyperboliques de dimension 3, *J. Diff. Geom.*, **35**, (1992), 609–626

[24] Labourie F., Problèmes de Monge-Ampère, courbes pseudo-holomorphes et laminations, *G.A.F.A*, **7**, (1997), 496–534

[25] Labourie F., Un Lemme de Morse pour les surfaces convexes, *Invent. Math.*, **141**, no. 2, (2000), 239–297

[26] Labourie F., Random k-surfaces, *Annals of Maths*, **161**, no. 1, (2005), 105-140

[27] Labourie F., Asymptotic counting of minimal surfaces in hyperbolic 3-manifolds, *Séminaire Bourbaki*, 05/2021, 73ème année, 2020–2021, no. 1179

[28] Omori H., Isometric immersions of riemannian manifolds, *J. Math. Soc. Japan*, **19**, no. 2, (1967), 205–214

[29] Petersen P., *Riemannian Geometry*, Graduate Texts in Mathematics, **171**, Springer Verlag, New York, (1998)

[30] Rosenberg H., Spruck J., On the existence of convex hypersurfaces of constant Gauss curvature in hyperbolic space, *J. Diff. Geom.*, **40**, no. 2, 379–409, (1994)

[31] Schlenker J. M., Hyperbolic manifolds with convex boundary, *Invent. Math.*, **163**, (2006), 109–169

[32] Smith G., Ph.D., Univ. Paris XI, Orsay, (2004)

[33] Smith G., An Arzela-Ascoli Theorem for Immersed Submanifolds, *Ann. Fac. Sci. Toulouse Math.*, **16**, no. 4, (2007), 817–866

[34] Smith G., Moduli of Flat Conformal Structures of Hyperbolic Type, *Geom. Dedicata*, **154**, no. 1, (2011), 47–80

[35] Smith G., Special Lagrangian curvature, *Math. Ann.*, **335**, no. 1, (2013), 57–95

[36] Smith G., The Plateau problem for convex curvature functions, *Ann. Inst. Fourier*, **70**, no. 1, (2020), 1–66

[37] Smith G., Möbius structures, hyperbolic ends and k-surfaces in hyperbolic space, arXiv:2104.03181